\documentclass[12pt,reqno]{amsart}
\usepackage{amsmath, amsthm, amscd, amssymb, amsfonts, latexsym}
\usepackage{fullpage}
\usepackage{bm}
\usepackage[all]{xy}
\usepackage{tikz}
\usepackage{mathrsfs}
\usepackage{dsfont}

\newcommand{\kommentar}[1]{}

\newcommand{\eps}{\varepsilon}

\newcommand{\Z}{\mathbb Z}
\newcommand{\Q}{\mathbb Q}

\newcommand{\co}{\mathcal O}
\newcommand{\C}{\mathbb C}

\newcommand{\sumstar}{\sideset{}{^*}\sum}

\newcommand{\ndv}{\nmid}
\newcommand{\edv}{\mathrel\Vert}

\DeclareMathOperator{\Gal}{Gal}

\DeclareMathOperator{\ord}{ord}

\renewcommand{\pmod}[1]{\,(\mathrm{mod}\,#1)}

\newtheorem{lem}{Lemma}[section]
\newtheorem{prop}[lem]{Proposition}
\newtheorem{thm}[lem]{Theorem}
\newtheorem{defn}[lem]{Definition}
\newtheorem{cor}[lem]{Corollary}

\theoremstyle{definition}

\newtheorem{rem}[lem]{Remark}

\author{Hung M. Bui, Alexandra Florea and Hieu T. Ngo}
\address{Department of Mathematics, University of Manchester, Manchester M13 9PL, UK}
\email{hung.bui@manchester.ac.uk}
\address{UC Irvine, Mathematics Department, Rowland Hall, Irvine 92697, USA}
\email{floreaa@uci.edu}
\address{Institute of Mathematics, Vietnam Academy of Science and Technology, Hanoi, Vietnam}
\email{nthieu@math.ac.vn}

\title{Simultaneous nonvanishing of Dirichlet $L$-functions in Galois orbits}

\allowdisplaybreaks

\begin{document}

\begin{abstract}
Under the Generalized Riemann Hypothesis, we prove that given any two  distinct imprimitive Dirichlet characters $\eta_1, \eta_2$ modulo $q=p^k$, a positive proportion of characters $\chi$ modulo $q$ in a fixed Galois orbit of primitive characters satisfies the nonvanishing property that $L(1/2,\chi \eta_1) L(1/2,\chi \eta_2) \neq 0$, as $k \to \infty$ (with $p$ fixed). Previously, only a positive proportion of nonvanishing result was available in Galois orbits (as opposed to simultaneously nonvanishing), due to work of Khan, Mili\'cevi\'c and Ngo. The main ingredients are obtaining a sharp upper bound on the mollified fourth moment over the Galois orbit using an Euler product mollifier, and obtaining a lower bound for the mollified second moment, which relies on using results from Diophantine approximation (such as the $p$-adic Roth theorem). We also unconditionally compute the second moments for $L$--functions associated to primitive Dirichlet characters in full orbits and thinner orbits. 
\end{abstract}

\maketitle

\section{Introduction and statement of results}

The values of $L$--functions at special points are central objects in number theory. It is expected that an $L$--function vanishes at the center of the critical strip either for deep, arithmetic reasons (such as the $L$--function associated to an elliptic curve, whose vanishing at the central point is related, through the Birch and Swinnerton-Dyer conjecture, to the rank of the elliptic curve) or for `obvious' reasons (such as when the root number appearing in the functional equation of the $L$--function is equal to $-1$). 

Chowla \cite{Cho65} predicted that $L(1/2,\chi) \neq 0$ for any Dirichlet character $\chi$. Much work has been done towards Chowla's conjecture. For the full family of Dirichlet $L$--functions, Iwaniec and Sarnak \cite{IS99} established that for at least $1/3-\varepsilon$ of the primitive Dirichlet characters $\chi$ modulo $q$, where $q$ is large enough in terms of $\varepsilon$, the values $L(1/2,\chi)$ are non-zero. To do that they computed mollified moments, using a one-piece mollifier. The proportion was recently improved to $7/19$ by Qin and Wu \cite{QW25}, using a two-piece mollifier (see also \cite{Bui12}). When $q$ is prime, the proportion is shown to be at least $5/13$ by Khan, Mili\'cevi\'c and Ngo \cite{KMN22} (see also \cite{KN16}). 

Restricting to the subset of quadratic characters, Soundararajan \cite{Sou00} computed mollified moments and proved that more than $7/8$ of $L(1/2,\chi)$ are non-zero, as $\chi$ ranges over primitive real characters. This proportion can be improved to $15/16$ if one assumes the Generalized Riemann Hypothesis (GRH) and computes the one-level density of zeros, as in the work of \"Ozl\"uk and Snyder \cite{OS99}. Higher order characters are more difficult to study, and the results are not as strong as for the family of quadratic characters. Focusing on the family of $L$--functions associated to cubic characters over $\mathbb{Q}$, Baier and Young \cite{BY10} showed that the number of primitive Dirichlet characters $\chi$ of order $3$ with conductor less than $Q$ for which $L(1/2,\chi) \neq 0$ is bounded below by $Q^{6/7-\varepsilon}$. Note that while this result provides an infinite family of $L$--functions that are nonvanishing at the central point, it does not imply a positive proportion of nonvanishing, since the family has size $Q$. Under GRH, David and G\"ulo\u{g}lu \cite{DG22} showed that more than $2/13$ of the $L$--functions of the  cubic Dirichlet characters $\chi_c$, for square-free $c$ in the Eisenstein field, do not vanish at the central point. Recently,  David, de Faveri, Dunn and Stucky \cite{DDDS24} removed the assumption of GRH and showed that more than $14\%$ of the central values in the same family of cubic $L$--functions are non-zero.

It is interesting to study the simultaneous question of nonvanishing of Dirichlet $L$--functions at the central point; there are fewer results in this direction in the literature. In \cite{MV02}, Michel and Vanderkam showed that given three Dirichlet characters $\chi_1,\chi_2,\chi_3$ of moduli $D_1,D_2, D_3$, a positive proportion of holomorphic Hecke cusp forms of weight $2$, prime level $q$ and trivial nebentypus are such that the central values $L(1,f \otimes \chi_j) \neq 0$, for $j=1,2,3$, for sufficiently large $q$. Focusing on the family of Dirichlet $L$--functions, Zacharias \cite{Zac19} showed that given two Dirichlet characters $\chi_1,\chi_2$ of prime modulus $q$, a positive proportion of $\chi$ modulo $q$ have the property that $L(1/2,\chi), L(1/2,\chi \chi_1), L(1/2,\chi \chi_2)$ are non-zero, for $q$ large enough.  

Instead of varying a Dirichlet character over all the primitive characters of the same modulus, there is a natural, more refined family of Galois orbit characters. More precisely, for $n$ a positive integer, let $\boldsymbol{\mu}_n$ denote the group of $n^{\rm th}$ roots of unity in $\overline{\mathbb{Q}}$. The Galois group $\text{Gal}(\mathbb{Q}(\boldsymbol{\mu}_{\phi(q)})/\mathbb{Q})$ acts naturally on the set of primitive Dirichlet characters modulo $q$ as follows. For $\sigma \in \text{Gal}(\mathbb{Q}(\boldsymbol{\mu}_{\phi(q)})/\mathbb{Q})$ and $\chi $ a primitive Dirichlet character modulo $q$, let $\chi^{\sigma}$ be the character defined by $\chi^{\sigma}(n) = \sigma(\chi(n))$. This Galois action partitions the set of primitive Dirichlet characters into orbits; each such Galois orbit is denoted by $\co$.

The algebro-analytic interest of the Galois conjugacy variation of $L$--values has roots in the works of Shimura \cite{Shi76,Shi77} and Rohrlich \cite{Roh80a,Roh80b}. 
Let $f(z)=\sum_{n=1}^\infty a_ne(nz)$ be a cusp form of even weight $2k$ and level $N$. For any automorphism $\sigma$ of $\C$, define the cusp form $f^\sigma$ by $f^\sigma(z)=\sum_{n=1}^\infty \sigma(a_n) e(nz)$. Then \cite[Theorem 1]{Shi77} implies that the central value $L(k,f\otimes\chi)\neq 0$ if and only if $L(k,f^\sigma\otimes\chi^\sigma)\neq 0$. Therefore, in order to single out a nonvanishing value $L(k,f\otimes\chi)\neq 0$, one can instead vary $\sigma$ in a finite subgroup $G$ of automorphisms of $\C$ and prove that an appropriate average of the values $\{L(k,f^\sigma\otimes\chi^\sigma): \sigma\in G\}$ is nonzero. In \cite{Roh80a,Roh80b}, the $L$-functions associated to certain classes of grossencharacters of an imaginary quadratic field were considered, and the appropriate average was taken to be an $L$-function associated to an ideal class in a ray class group (see also \cite{Roh80c,MR82}). Later in the 1980s, strong results concerning the order of vanishing of central $L$--values were established and deep consequences were then derived for the arithmetic of rational (i.e.~defined over $\mathbb{Q}$) elliptic curves, more particularly those with complex multiplication. We now give a brief sketch of this line of inquiry.

Let $E$ be a rational elliptic curve of conductor $N_E$, say, with complex multiplication by the ring of integer $\co_K$ in an imaginary quadratic field $K$. To $E$ there is an associated grossencharacter $\psi_E$ such that the Hecke $L$--function
$$
L(s,\psi_E)=\sum_{\mathfrak{a}\subset\co_K} \frac{\psi_E(\mathfrak{a})}{\mathcal{N}_{K/\Q}(\mathfrak{a})^{s}}
$$
can be identified with the Hasse-Weil $L$--function of $E/\mathbb{Q}$, i.e.~$L(s,E/\mathbb{Q})=L(s,\psi_E)$. This $L$-function has analytic continuation and the functional equation
\begin{align*}
\Lambda(s,E/\mathbb{Q})&:=\bigg(\frac{N_E}{4\pi^2}\bigg)^{s/2}\Gamma(s)L(s,E/\mathbb{Q})\\
&=w_{E}\Lambda(2-s,E/\mathbb{Q})
\end{align*}
with the root number $w_{E}=\pm1$. 

For any positive integer $k$, the $L$--function $L(s,\psi_E^{2k+1})$ satisfies the functional equation
\begin{align*}
\Lambda(s,\psi_E^{2k+1})&:=\Bigg(\frac{N_{E}^{(2k+1)}}{4\pi^2}\Bigg)^{s/2}\Gamma(s)L(s,\psi_E^{2k+1})\\
&=w_{E}^{(2k+1)}\Lambda(2k+2-s,\psi_E^{2k+1})
\end{align*}
with $w_{E}^{(2k+1)}=\pm1$. Here $N_{E}^{(2k+1)}=|{\rm disc}(K)|\mathcal{N}_{K/\Q}(f_{\psi_E^{2k+1}})$ where $f_{\psi_E^{2k+1}}$ is the conductor of $\psi_E^{2k+1}$. 
Greenberg \cite[Theorem 4]{Gre83} showed that the trivial critical zeros account for all but finitely many vanishing central values $L(k+1,\psi_E^{2k+1})$. In \cite{Gre85}, the same conclusion was derived for the vanishing of the central values $L(k+1,\psi_{2k+1})$ where $\psi_{2k+1}$ varies over all grossencharacters of an imaginary quadratic field $K$ of type $A_0$ with infinity type $(2k+1,0)$ and whose conductors have bounded norm.
In \cite{Roh84b}, dropping the complex multiplication assumption for $E$ and considering a finite set $S$ of primes not dividing $N_E$, Rohrlich showed that when $\chi$ ranges over all primitive Dirichlet characters whose conductors are divisible by only primes in $S$, the central values $L(1,E\otimes\chi)$ are nonvanishing for all but finitely many such $\chi$. For this, Rohrlich varied $\chi$ in a Galois orbit and computed the first moment of the central values of the family $\{L(1,E\otimes\chi):\chi\in\co\}$. Chinta \cite{Chi02} computed the mollified first moment for the same family as in \cite{Roh84b}, extending to those $\chi$ of prime moduli (see \cite{KMN16} for a comparison of their methods). 

In \cite{Roh84a}, Rohrlich considered twists of $\psi_E$ by characters of anti-cyclotomic extensions of $K$. More precisely, let $S$ be a finite set of primes and let $K_{S,\infty}^-$ be the compositum of all the anti-cyclotomic $\Z_l$-extensions $K_{l,\infty}^-$ of $K$, where $l$ ranges over $S$. It turns out that ${\rm Gal}(K_{S,\infty}^-/K)$ is isomorphic to the product of a finite group and $\prod_{l\in S} \Z_l$. A (finite order) character $\rho$ of ${\rm Gal}(K_{S,\infty}^-/K)$ can be identified with a Hecker character (i.e.~an idele class character) of $K$; let $\rho\psi_E$ denote the primitive Hecker characters determined by the product of $\rho$ and $\psi_E$. Let $X$ denote the set of Hecke characters $\chi=\rho\psi_E$ obtained in this way. In simpler terms, one can let $\rho$ vary over all grossencharacters of finite order of $K$ whose conductors are divisible only by primes above those in $S$ such that $\rho\circ c=\rho^{-1}$. Here $c$ denotes the complex conjugation in $K$. 
If the root number $w_{\chi}=1$, then $L(1,\chi)\neq 0$ for all but finitely many such $\chi\in X$. If $w_{\chi}=-1$, then $L(s,\chi)$ has a simple zero at $s=1$ for all but finitely many such $\chi\in X$. 

Amazingly, Roth's theorem \cite{Rot55} was used in \cite{Gre85} and its non-archimedean analog (proved by Ridout \cite{Rid58}) occurred in \cite{Roh84a} as important ingredients. We shall refer to the $p$-adic Roth's theorem in \cite{Rid58} as the Roth-Ridout theorem.

The Galois conjugacy variation and the use of deep results in Diophantine approximation in computing moments of $L$--values not only possess inherent interest, but
also find important arithmetic implications for the Mordell-Weil groups of rational points of $E$ with values in towers of number fields, for instance $\Z_l$-extensions of $K$. Let $K_{l,\infty}^{-}=\bigcup_{n\geq 0}K_{l,n}^-$ be the anti-cyclotomic $\mathbb{Z}_l$-extension of $K$ with the $n^{\rm th}$-layer $K_{l,n}^{-}$ being a $\mathbb{Z}/l^n\mathbb{Z}$-extension of $K$ and a dihedral Galois extension over $\Q$. 
Let $K_{l,\infty}=\bigcup_{n\geq 0}K_{l,n}$ be any $\mathbb{Z}_l$-extension of $K$ with the $n^{\rm th}$-layer $K_{l,n}$.
Assume that $l>3$; the torsion subgroup of $E(K_{l,\infty})$ can be described explicitly \cite{Avc25}, which is finite in particular.
Greenberg \cite{Gre83,Gre87,Gre01} showed that if $K_{l,\infty}\neq K_{l,\infty}^{-}$, then the $\Z$-ranks of $E(K_{l,n})$ stabilize as $n\to\infty$, whence $E(K_{l,\infty})$ is finitely generated. 
Assume also that $E$ has good ordinary reduction at $l$ and consider the anti-cyclotomic extension $K_{l,\infty}^{-}/K$. 
If $w_{E}=1$, then $E(K_{l,\infty}^{-})$ is finitely generated \cite{Gre83,Gre01}.
If $w_{E}=-1$, then the $\Z$-ranks of $E(K_{l,n}^{-})$ are unbounded and equal $2p^n+c$, where $c$ is a constant, for all sufficiently large $n$ \cite{Gre01}. With $S$ a finite set of primes where $E$ has good reduction, let $\Q_S^{\rm ab}$ be the maximal abelian extension of $\Q$ unramified outside $S$ and infinity; in other words $\Q_S^{\rm ab}$ is generated over $\Q$ by the roots of unity of $l$-power order, where $l$ ranges over $S$. Rohrlich \cite{Roh84b} showed that $E(\Q_S^{\rm ab})$ is finitely generated. Those results were all directed towards the question raised by Mazur in \cite{Maz72} (see \cite{Roh88} for a generalization in this direction).


In this paper, we are interested in the question of the simultaneous nonvanishing of Dirichlet $L$--functions of prime power moduli in Galois orbits at the central point. We note that for every divisor $c|(p-1)$, there is a unique orbit $\co=\co_c$ consisting of all primitive characters $\chi$ modulo $q = p^k$ of order $p^{k-1}(p-1)/c$, and having size $\phi(p^{k-1}(p-1)/c)$. 
In this context, Khan, Mili\'cevi\'c and Ngo \cite{KMN16} proved that at least $1/3$ of the $L$--functions in each orbit are non-zero at the central point, as $q=p^k \to \infty$ (with $p$ fixed), recovering the nonvanishing proportion obtained by Iwaniec and Sarnak \cite{IS99} for the whole family of Dirichlet characters. More precisely, it was proved in \cite{KMN16} that if $q=p^k$, then
\[
\frac{1}{|\co|}\sum_{\substack{\chi\in\co\\L(1/2,\chi)\ne0}}1\geq \frac13-o(1)
\]
as $p$ is fixed and $k\rightarrow\infty$. In the other direction, Munsch and Shparlinski \cite{MM23} studied an analogous problem when $k=1$ and $p\to\infty$, and obtained some nonvanishing results, conditional on the size of solutions of certain linear congruences being large enough.

Let $\eta_1$, $\eta_2$ be two distinct imprimitive Dirichlet characters modulo $q$. We are interested in the simultaneous nonvanishing of $L(1/2,\chi \eta_1)$ and $L(1/2,\chi \eta_2)$, as $\chi$ ranges over a Galois orbit $\co$ of primitive characters modulo $q$, and as $k \to \infty$.  
The main result that we prove in this paper is the following nonvanishing result.

\begin{thm}\label{thm:nonvanishing-full-orbit}
Assume GRH. Let $q=p^k$ for a fixed odd prime $p$. Let $\co$ be any Galois orbit of primitive Dirichlet characters modulo $q$, and let $\eta_1, \eta_2$ be two  distinct imprimitive Dirichlet characters modulo $q$.  We have
\begin{align*}
&\#\big\{\chi\in\co:L(\tfrac12,\chi\eta_1)L(\tfrac12,\chi\eta_2)\ne0\big\}\gg q.
\end{align*}
\end{thm}

Proving the theorem above relies on estimating mollified moments of $L$--functions; namely, obtaining a lower bound for a mollified second moment, and obtaining a sharp upper bound for a mollified fourth moment. The choice of mollifier is a key step in obtaining a sharp upper bound for the mollified fourth moment. Namely, we construct an Euler product mollifier, rather than a classical Dirichlet series mollifier. We note that similar mollifier constructions have been recently used in many settings to prove a variety of results. Lester and Radziwi\l\l \,  introduced the Euler product mollifier in \cite{LR21} to study sign changes of Fourier coefficients of half-integral weight modular forms. This construction was subsequently used to obtain nonvanishing results for $L$--functions \cite{{DFL21}} and to obtain weighted central limit theorems for central values of $L$--functions \cite{BELP25}, among many other results. 

Let $M(\chi)$ be a mollifier as defined in Section \ref{mol_section} (see equations \eqref{mollifier} and \eqref{mj}; also \eqref{beta_j}, \eqref{condition_c}, \eqref{sj}). Proving Theorem \ref{thm:nonvanishing-full-orbit} relies on the following two key estimates. The first estimate provides a lower bound for the mollified second moment (for the purpose of Theorem \ref{mol_lb} below, in the definition of the mollifier, we take $v=4$). 

\begin{thm}
\label{mol_lb}
Under the same assumptions as in Theorem \ref{thm:nonvanishing-full-orbit}, we have that
$$ \frac{1}{|\co|}\,\sum_{\substack{\chi\in\co}}L(\tfrac12,\chi\eta_1)L(\tfrac12,\overline{\chi\eta_2})M(\chi\eta_1)M(\overline{\chi\eta_2}) \gg 1.$$
\end{thm}

We also have the following upper bound.

\begin{thm}
\label{mol_ub}
Assume GRH. For any integer $v \geq 1$, we have that 
$$\sumstar_{\chi \pmod q} \big| L (\tfrac{1}{2},\chi) M(\chi) \big|^v \ll q,$$
where $\sumstar$ denotes summation over primitive characters.
\end{thm}

\begin{rem}
We note that the upper bound in Theorem \ref{mol_ub} holds for all positive integers $q \to \infty$, not neccesarily prime powers as in Theorems \ref{thm:nonvanishing-full-orbit} and \ref{mol_lb}.
\end{rem}
The lower bound in Theorem \ref{mol_lb} uses the following twisted moment asymptotic. 

\begin{thm}\label{thm:twisted-moment-full-orbit}
Let $q=p^k$ for a fixed odd prime $p$ and $\co$ be any Galois orbit of primitive Dirichlet characters modulo $q$. 
Suppose that $\eta_1,\eta_2$ are two distinct imprimitive Dirichlet characters modulo $q$; let $\eta_1\overline{\eta_2}$ be induced from a primitive character modulo $p^h$ with $0<h<k$, say.
Let $\vartheta\geq 0$ and $m_1,m_2$ be two positive integers; put $m_1'=m_1/(m_1,m_2)$ and $m_2'=m_2/(m_1,m_2)$. We have that
\begin{align*}
&\frac{1}{|\co|}\sum_{\substack{\chi\in\co}}L(\tfrac12,\chi\eta_1)L(\tfrac12,\overline{\chi\eta_2})\chi(m_1)\overline{\chi}(m_2)=\frac{\eta_1(m_2')\overline{\eta_2}(m_1')}{\sqrt{m_1'm_2'}}L(1,\eta_1\overline{\eta_2})+E(m_1,m_2) ,
\end{align*}
where $E(m_1,m_2)$ satisfies the following property: 
for any $\varepsilon>0$ and any sequence of complex numbers $(x_m)_{m=1}^\infty$ satisfying $x_m\ll_\varepsilon m^\varepsilon$, if $k$ is large enough in terms of $p$ and $\varepsilon$ then
$$
\mathcal{E}(\vartheta):=\sum_{m_1,m_2\leq q^\vartheta}\frac{x_{m_1}\overline{x_{m_2}}}{\sqrt{m_1m_2}}E(m_1,m_2) 
    \ll_\varepsilon p^{h/4} q^{-1/2+\vartheta+\varepsilon} + p^{-h/4}q^{-1/4+\varepsilon}.
$$
\end{thm}

We have as consequence a second moment asymptotic along a full Galois orbit.

\begin{cor}\label{cor:full-orbit-moment}
Under the same assumptions as in Theorem \ref{thm:twisted-moment-full-orbit}, we have
\begin{align*}
&\frac{1}{|\co|}\sum_{\substack{\chi\in\co}}L(\tfrac12,\chi\eta_1)L(\tfrac12,\overline{\chi\eta_2})=L(1,\eta_1\overline{\eta_2})+
O_\varepsilon(p^{h/4} q^{-1/2+\varepsilon} + p^{-h/4}q^{-1/4+\varepsilon}).
\end{align*}
\end{cor}

\begin{rem}
We note that Theorem \ref{thm:twisted-moment-full-orbit} and Corollary \ref{cor:full-orbit-moment} are unconditional, while Theorem \ref{mol_lb} is conditional on GRH. This is due to the fact that, in order to obtain a lower bound on the mollified moment, one needs a strong lower bound on the size of $L(1,\eta_1 \overline{\eta_2})$, which is currently beyond reach using the classical zero-free region of $L(s,\eta_1 \overline{\eta_2})$. It is possible that one could obtain Theorem \ref{thm:twisted-moment-full-orbit} without the full force of GRH, but rather under the assumption of a good zero-free region. For simplicity (and to avoid the existence of possible Landau-Siegel zeros), we have assumed GRH.
\end{rem}
We can also derive an unconditional (weaker) simultaneous nonvanishing result.

\begin{cor}\label{cor:full-orbit-nonvanishing}
Under the same assumptions as in Theorem \ref{thm:twisted-moment-full-orbit}, we have
\begin{align*}
&\#\big\{\chi\in\co:L(\tfrac12,\chi\eta_1)L(\tfrac12,\chi\eta_2)\ne0\big\}\gg_\varepsilon q^{2/3-\varepsilon}.
\end{align*}
\end{cor}

We now formulate analogous problem for the so-called thin Galois orbits.
For $0\leq h\leq k$, we put $q_h=\phi(p^h)=p^{h-1}(p-1)$.
Inside the number field $\mathbb{Q}(\boldsymbol{\mu}_{q_k})$, there is a natural tower of number fields
$$
\Q(\boldsymbol{\mu}_{q_k}) := K_{k-1} \supsetneq \cdots \supsetneq K_0 \supsetneq \Q.
$$
where $K_{k-1-\kappa}:=\Q(\boldsymbol{\mu}_{q_{k-\kappa}})$ for $0\leq \kappa \leq k-1$. 
Let $\co$ be a Galois orbit of primitive Dirichlet characters modulo $q$ under the action of $\text{Gal}(\mathbb{Q}(\boldsymbol{\mu}_{q_k})/\mathbb{Q})$. 
The subgroup 
$\text{Gal}(\mathbb{Q}(\boldsymbol{\mu}_{q_k})/K_{k-1-\kappa})$
of $\text{Gal}(\mathbb{Q}(\boldsymbol{\mu}_{q_k})/\mathbb{Q})$ acts on $\co$; we call an orbit under this action a \emph{thin Galois orbit} and denote by $\co_\kappa$ any such orbit. For every $0\leq \kappa \leq k-1$, all thin Galois orbits $\co_\kappa$ are of equal size
$$
|\co_\kappa|= 
    \begin{cases}
    p^\kappa & \text{ if }\, 0\leq \kappa < k-1, \\
    q_{k-1} & \text{ if }\, \kappa = k-1 .
    \end{cases}
$$

We have an unconditional analog of Theorem \ref{thm:twisted-moment-full-orbit} for thin Galois orbits.

\begin{thm}\label{thm:twisted-moment-thin-orbit}
Let $q=p^k$ for a fixed odd prime $p$, $k-h/2<\kappa\leq k-1$, and $\co_\kappa$ be any Galois orbit of primitive Dirichlet characters modulo $q$ under the action of $\text{Gal}(K_{k-1}/K_{k-1-\kappa})$. 
Let $\eta_1, \eta_2$ be two distinct imprimitive Dirichlet characters modulo $q$. Suppose that $\eta_1\overline{\eta_2}$ is induced from a primitive character modulo $p^h$ with $0<h<k$.
Let $\vartheta\geq 0$ and $m_1,m_2$ be two positive integers; put $m_1'=m_1/(m_1,m_2)$ and $m_2'=m_2/(m_1,m_2)$. 
We have that
\begin{align*}
&\frac{1}{|\co_\kappa|}\sum_{\substack{\chi\in\co_\kappa}}L(\tfrac12,\chi\eta_1)L(\tfrac12,\overline{\chi\eta_2})\chi(m_1)\overline{\chi}(m_2)=\frac{\eta_1(m_2')\overline{\eta_2}(m_1')}{\sqrt{m_1'm_2'}}L(1,\eta_1\overline{\eta_2})+E_\kappa(m_1,m_2) ,
\end{align*}
where $E_\kappa(m_1,m_2)$ satisfies the following property: for any $\varepsilon>0$ and any sequence of complex numbers $(x_m)_{m=1}^\infty$ satisfying $x_m\ll_\varepsilon m^\varepsilon$, if $k$ is large enough in terms of $p$ and $\varepsilon$, then
$$
\mathcal{E}_\kappa(\vartheta):=\sum_{m_1,m_2\leq q^\vartheta}\frac{x_{m_1}\overline{x_{m_2}}}{\sqrt{m_1m_2}}E_\kappa(m_1,m_2)
\ll_\varepsilon 
p^{h/4-\kappa} q^{1/2+\vartheta+\varepsilon} + p^{-\kappa/4}q^\varepsilon 
+ p^{-\kappa+(3k-h)/4}q^{\varepsilon}.
$$
\end{thm}

The error term estimate in Theorem \ref{thm:twisted-moment-thin-orbit} can be made admissible, say 
$\mathcal{E}_\kappa(\vartheta) = O_{\varepsilon}(q^{-\varepsilon})$,
by assuming that $(k,h,\kappa)$ is $(\vartheta,\varepsilon)$-compatible, where the compatibility is defined as follows.

\begin{defn}
Let $k,h,\kappa$ be positive integers. For $\vartheta\geq0$ and $\varepsilon>0$, we say that the datum $(k,h,\kappa)$ is \emph{$(\vartheta,\varepsilon)$-compatible} if
\begin{enumerate}
    \item $k-\frac{h}{2}+\varepsilon h<\kappa<k$,  
    \item $\kappa> \frac{h}{4} + (1/2+\vartheta+\varepsilon)k$.
\end{enumerate}
\end{defn}

We have an unconditional second moment asymptotic along a thin Galois orbit.

\begin{cor}\label{cor:thin-orbit-moment}
Under the same assumptions as in Theorem \ref{thm:twisted-moment-thin-orbit}, for any $\varepsilon>0$, if $(k,h,\kappa)$ is $(0,\varepsilon)$-compatible then 
\begin{align*}
&\frac{1}{|\co_\kappa|}\sum_{\substack{\chi\in\co_\kappa}}L(\tfrac12,\chi\eta_1)L(\tfrac12,\overline{\chi\eta_2})=L(1,\eta_1\overline{\eta_2})+O_\varepsilon(q^{-\varepsilon}).
\end{align*}
\end{cor}

\section{Preliminaries}

\subsection{Notation}

    Throughout the paper, unless otherwise stated, we let $p$ denote a fixed odd prime and $q=p^k$ with $k\geq 1$. 

    A Dirichlet character $\chi$ modulo $q$ is said to have \emph{height} $0\leq h\leq k$ if it is induced from a primitive Dirichlet character $\chi^*$ modulo $p^h$.
    For $j\in\{1,2\}$, let $\eta_j$ be an imprimitive Dirichlet character modulo $q$ of height $0\leq  h_i < k$, and let $\kappa_j\in\{0,1\}$ be such that $\chi\eta_j(-1)=(-1)^{\kappa_j}$.
    Write $\eta=\eta_1\overline{\eta_2}$, and denote by $0\leq h < k$ the height of $\eta$; in particular $0\leq h\leq \max(h_1,h_2) < k$. 
    If $\eta_1\neq\eta_2$, then $h\geq 1$.

    When summing over integers, $\sumstar$ indicates that $p$ does not divide the summation variable(s). We also use the notation $\sumstar$ when summing over primitive characters $\chi$ modulo $q$. The notations shall not be confused, since they are over different objects. We write $\ord_p(n)$ for the biggest integer $\nu$ such that $p^\nu | n$. The \emph{multiplicative order} $\ord(n;q)$ is the order of $n$ in the multiplicative group $(\Z/q\Z)^\times$.
    The Euler totient function is denoted by $\phi(\cdot)$. For $0\leq h\leq k$, we put $q_h=\phi(p^h)=p^{h-1}(p-1)$. We write $\Omega(n)$ for the number of prime divisors of $n$ counted with multiplicity. The Liouville function is $\lambda(\prod_p p^{a_p})=(-1)^{\sum_p a_p}$.
    
    We use $\co$ to denote a \emph{full Galois orbit}, i.e.~an orbit of primitive Dirichlet characters modulo $q$ under the action of $\Gal\big(\Q(\boldsymbol{\mu}_{\phi(q)})/\Q\big)$. For $0\leq \kappa \leq k-1$, we write $\co_\kappa$ for an orbit of primitive Dirichlet characters modulo $q$ under the action of $\Gal\big(\Q(\boldsymbol{\mu}_{\phi(q)})/\Q(\boldsymbol{\mu}_{\phi(p^{k-\kappa})})\big)$, and call it a \emph{thin Galois orbit}.
    We make use of the expectation notation for the character average over an orbit (full or thin), namely
    \begin{equation}\label{eq:expectation-orbit-character}
        \mathbb{E}_\co[\chi(n)] = \frac{1}{|\co|}  \sum_{\chi\in \co}  \chi(n), \quad
        \mathbb{E}_{\co_\kappa}[\chi(n)] = \frac{1}{|\co_\kappa|} \sum_{\chi\in \co_\kappa} \chi(n).     
    \end{equation}

    For natural numbers $n,r$ and $\zeta\in\boldsymbol{\mu}_{p-1}$, we write the congruence $n\equiv \zeta \ (\text{mod}\ p^r)$ to indicate that $n-\zeta\in p^r\Z_p$.

We define the normalized Gauss sum $\epsilon_\chi$ associated to a Dirichlet character $\chi$ modulo $q$ by 
\begin{equation}
\label{normalized}
\epsilon_\chi=\frac{c_\chi(1)}{\sqrt{q}},
\end{equation}
where the generalized Gauss sum $c_\chi$ is given by
$$
c_\chi(n) = \sumstar_{u (\text{mod}\ q)}\chi(u)e_q(un)
$$
and $e_q(x)=e(\frac{x}{q})$.

    We write $\varepsilon$ for a small positive constant and use the convention that one $\varepsilon$ may differ from another $\varepsilon$, but only within an absolute constant multiple. More precisely, two $\varepsilon$ notations are indeed two positive constants $\varepsilon_1$ and $\varepsilon_2$ such that there are absolute constants $c'>c>0$ with $c\varepsilon_1<\varepsilon_2<c'\varepsilon_1$.

\subsection{Approximate functional equation}

Let $\chi$ be a primitive Dirichlet character modulo $q$. Our assumptions imply that $\chi\eta_j$ is primitive modulo $q$ for $j\in\{1,2\}$, and so the Dirichlet $L$-function $L(s,\chi\eta_j)$ satisfies the functional equation
\begin{align}\label{fe}
\Lambda(\tfrac12+s,\chi\eta_j)&:=\Big(\frac{q}{\pi}\Big)^{s/2}\Gamma\Big(\frac{1/2+s+\kappa_j}{2}\Big)L(\tfrac12+s,\chi\eta_j)\nonumber\\
&=i^{-\kappa_j} \epsilon_{\chi\eta_j}\Lambda(\tfrac12-s,\overline{\chi\eta_j}).
\end{align}

We are in a position to derive a standard approximate functional equation in our context.

\begin{lem}\label{lem:afe}
Let $G(s)$ be an even entire function of rapid decay in any fixed strip $|\emph{Re}(s)| \leq C$ satisfying $G(0)= 1$.
For $x>0$, define
\begin{equation}\label{formulaV+}
V(x)=\frac{1}{2\pi i}\int_{(1)}\frac{\Gamma\big(\frac{1/2+s+\kappa_1}{2}\big)\Gamma\big(\frac{1/2+s+\kappa_2}{2}\big)}{\Gamma\big(\frac{1/2+\kappa_1}{2}\big)\Gamma\big(\frac{1/2+\kappa_2}{2}\big)}G(s)x^{-s}\frac{ds}{s}.
\end{equation}
For any $X>0$, we have
\begin{align*}
L(\tfrac 12,\chi\eta_1)L(\tfrac 12,\overline{\chi\eta_2})&=\sum_{m,n\geq1}\frac{\chi\eta_1(m)\overline{\chi\eta_2}(n)}{\sqrt{mn}}V\Big(\frac{\pi mn}{qX}\Big)\\
&\qquad + i^{-(\kappa_1+\kappa_2)}\epsilon_{\chi\eta_1}\epsilon_{\overline{\chi\eta_2}}\sum_{m,n\geq1}\frac{\overline{\chi\eta_1}(m)\chi\eta_2(n)}{\sqrt{mn}}V\Big(\frac{\pi mnX}{q}\Big).\nonumber
\end{align*}
\end{lem}

\begin{proof}
From Cauchy's residue theorem we have
\begin{align}\label{residue}
&\frac{1}{2\pi i}\int_{(1)}G(s)X^s\frac{\Lambda(1/2+s,\chi\eta_1)\Lambda(1/2+s,\overline{\chi\eta_2})}{\Gamma(\frac{1/2+\kappa_1}{2})\Gamma(\frac{1/2+\kappa_2}{2})}\frac{ds}{s}\nonumber\\
&\qquad=\frac{1}{2\pi i}\int_{(-1)}G(s)X^s\frac{\Lambda(1/2+s,\chi\eta_1)\Lambda(1/2+s,\overline{\chi\eta_2})}{\Gamma(\frac{1/2+\kappa_1}{2})\Gamma(\frac{1/2+\kappa_2}{2})}\frac{ds}{s}+\text{Res}_{s=0}.
\end{align}
The residue at $s=0$ is $L(\tfrac 12,\chi\eta_1)L(\tfrac 12,\overline{\chi\eta_2})$.
By \eqref{fe} and the change of variables $s\longleftrightarrow-s$ in \eqref{residue}, we obtain 
\begin{align*}
&L(\tfrac 12,\chi\eta_1)L(\tfrac 12,\overline{\chi\eta_2})=\frac{1}{2\pi i}\int_{(1)}G(s)X^s\frac{\Lambda(1/2+s,\chi\eta_1)\Lambda(1/2+s,\overline{\chi\eta_2})}{\Gamma(\frac{1/2+\kappa_1}{2})\Gamma(\frac{1/2+\kappa_2}{2})}\frac{ds}{s}\\
&\qquad\ + i^{-(\kappa_1+\kappa_2)}\epsilon_{\chi\eta_1}\epsilon_{\overline{\chi\eta_2}}\frac{1}{2\pi i}\int_{(1)}G(s)X^{-s}\frac{\Lambda(1/2+s,\overline{\chi\eta_1})\Lambda(1/2+s,\chi\eta_2)}{\Gamma(\frac{1/2+\kappa_1}{2})\Gamma(\frac{1/2+\kappa_2}{2})}\frac{ds}{s}.
\end{align*}
The lemma now follows by writing the $L$-functions in terms of Dirichlet series and then integrating term-by-term.
\end{proof}

\begin{rem}\label{boundforV}
We have
\[
x^j \frac{\partial^jV(x)}{\partial x^j}\ll_{j,C}(1+x)^{-C}
\]
for any fixed $j\geq0$ and $C > 0$. Furthermore, we have
\[
V(x)=1+O_\varepsilon(x^{1/2-\varepsilon}).
\]
\end{rem}


\subsection{Character averages for full Galois orbits}

In this section, we describe full Galois orbits explicitly and study their character averages.

For $1\leq h\leq k$, the multiplicative group $D_h:=(\Z/p^h\Z)^\times$ 
is cyclic of order $q_h$. There is a short exact sequence of group homomorphisms
$$
1\longrightarrow N^{(h)}_k \longrightarrow D_k \longrightarrow D_h \longrightarrow 1
$$
where $N^{(h)}_k:=\{1+p^hr: r\ (\text{mod}\ p^{k-h})\}$.
Let $g_k$ be a generator of $D_k$ and choose a generator $g_h$ for $D_h$ such that $g_h$ is the image of $g_k$ under the projection map. For $n\in D_k$, write $\iota_k(n)=\text{ind}_{g_k}(n)$ for the index of $n$ with respect to $g_k$. An element $n \in D_k$ belongs to $N^{(h)}_k$ if and only if $q_h|\iota_k(n)$, i.e. $\iota_k(n)=q_hr$ for some $r\ (\text{mod}\ p^{k-h})$.

The choice of a generator $g_k$ induces an isomorphism between the character group $C_k:=\widehat{(\Z/q\Z)^\times}$ and the cyclic group $\Z/q_k\Z$ as follows. Each character $\chi\in C_k$ is of the form
$$\chi_a(n)=e\left(\frac{a\iota_k(n)}{q_k}\right) \quad (n\in D_k)$$
for a unique $a=a_\chi\ (\text{mod}\ q_k)$, called the \emph{residue modulo $q_k$} of $\chi$. Recall that $\chi$ has conductor $q_h$ for some $1\leq h\leq k$ if $\chi$ is induced from a primitive Dirichlet character modulo $p^h$; hence, $\chi$ has \emph{height} $h$. It follows that $\chi$ has height $h$ if and only if $h$ is the biggest integer such that the restriction of $\chi$ to $N^{(h)}_k$ is trivial, or equivalently $p^{k-h} \edv a_\chi$. In particular, $\chi$ is primitive if and only if $p\nmid a_\chi$. Since $a_{\chi\eta}=a_\chi+a_\eta$, we see that the product of a primitive character $\chi$ with an imprimitive character $\eta$ is primitive. Note that the product of two primitive characters is not primitive in general.

The Galois group $\text{Gal}(\mathbb{Q}(\boldsymbol{\mu}_{\phi(q)})/\mathbb{Q})$ acts naturally on the set of primitive Dirichlet characters modulo $q$ by acting on the values of the characters. This action partitions the set of primitive Dirichlet characters into orbits; each such orbit is denoted by $\co$.
If $\chi=\chi_a$ is primitive with $a$ the residue of $\chi$, we put $c_\chi=\gcd(a,q_k)=\gcd(a,p-1)$ and call $c_\chi$ the \emph{characteristic} of $\chi$. For every divisor $c|(p-1)$, there is a unique orbit $\co=\co_c$ consisting of all primitive characters $\chi$ with characteristic $c_\chi=c$. Each character in $\co$ has the same order $q_k/c$, and the orbit $\co$ has size $|\co|=\phi(q_k/c)$.

\begin{lem}\label{lem:full-orbit-char-avg}
Let $n\in (\Z/q\Z)^\times$. If $\mathbb{E}_\co[\chi(n)]\neq 0$, then $n^{p-1} \equiv 1 \ (\emph{mod}\ p^{k-1})$. 
\end{lem}

\begin{proof}
The orbit character $\mathbb{E}_\co[\chi(n)]$ defined in \eqref{eq:expectation-orbit-character} satisfies the identity (see \cite[p.17]{Chi02} or \cite[Equation (2.8)]{KMN16}) 
\begin{align*} 
\mathbb{E}_\co[\chi(n)] = \frac{\mu(\ord(n^c;q))}{\phi(\ord(n^c;q))}. 
\end{align*}
If
$\mathbb{E}_\co[\chi(n)]\neq 0$, then $p^2\ndv \ord(n^c;q)$, equivalently $p^{k-2}|\iota_k(n^c)$. Therefore, $\mathbb{E}_\co[\chi(n)]\neq 0$ yields $n^{p(p-1)} \equiv 1 \ (\text{mod}\ p^k)$, whence $n^{p-1} \equiv 1 \ (\text{mod}\ p^{k-1})$ by Lemma \ref{lem:congruences-prime-power}. 
In other words, $\mathbb{E}_\co[\chi(n)]\neq 0$ only if $n\equiv \zeta \ (\text{mod}\ p^{k-1})$ for some $\zeta\in\boldsymbol{\mu}_{p-1}$. The lemma is proved.
\end{proof}

The following lemma evaluates $c_\eta$ when $\eta$ is an imprimitive character.

\begin{lem}\label{lem:generalized_Gauss_sum}
Let $n\in \Z/q\Z$. 
Let $\nu={\rm ord}_p(n)$, and write $n=p^\nu n'$ with $\gcd(p,n')=1$. 
Suppose that $\eta$ is a nontrivial imprimitive character modulo $q$ of height $0 < h < k$; in other words $\eta$ is induced from a primitive character $\eta^*$ modulo $p^h$.
Then, $c_\eta(n)$ is nonzero if and only if $\nu=k-h$, in which case
    $$
    c_\eta(n)=p^{k-h} \tau(\eta^*) \overline{\eta^*}(n') \quad \text{and} \quad |c_\eta(n)|=p^{k-h/2}.
    $$
\end{lem}

\begin{proof}
    By \cite[Theorem 9.12]{MV07}, $c_\eta(n)$ is nonzero only if $\nu \leq k-h$. We now assume that $\nu \leq k-h$. Then by the same theorem,
    $$
    c_\eta(n)=\overline{\eta^*}(n') \eta^*(p^{k-h-\nu}) \mu(p^{k-h-\nu})
    \frac{\varphi(p^k)}{\varphi(p^{k-\nu})} \tau(\eta^*).
    $$
    The factor $\eta^*(p^{k-h-\nu})$ is nonzero if and only if $k-h-\nu=0$, in which case
    $$
    c_\eta(n)=p^{k-h} \tau(\eta^*) \overline{\eta^*}(n').
    $$
    It follows that $|c_\eta(n)|=p^{k-h+h/2}=p^{k-h/2}$.
    The lemma is proved.
\end{proof}

We now compute a character average that will appear in the moment calculation.

\begin{lem}\label{lem:averageing_character_twist_chieta1_barchieta2}
    Let $n\in (\Z/q\Z)^\times$.
            If $n^{p-1} \not\equiv 1 \ (\emph{mod}\ p^{k-h})$, then
    $$
    \frac{1}{|\co|}  \sum_{\chi\in \co} \epsilon_{\chi\eta_1}\epsilon_{\overline{\chi\eta_2}} \chi(n) = 0.
    $$
    If $n^{p-1} \equiv 1 \ (\emph{mod}\ p^{k-h})$, then
    $$
    \frac{1}{|\co|}  \sum_{\chi\in \co} \epsilon_{\chi\eta_1}\epsilon_{\overline{\chi\eta_2}} \chi(n) \ll p^{-h/2}.
    $$
\end{lem}

\begin{proof}
   Using equation \eqref{normalized} and opening both Gauss sums, we have
        \begin{align*}
    	\frac{1}{|\co|}  \sum_{\chi\in \co} \epsilon_{\chi\eta_1}\epsilon_{\overline{\chi\eta_2}} \chi(n)
    	&= \frac{1}{q|\co|} \sum_{\chi\in \co}
    	\chi(n)\ \sumstar_{u,v \, (\text{mod}\ q)} \chi\eta_1(u)\overline{\chi\eta_2}(v)e_q(u+v) \\
    	&=  \frac{1}{q|\co|} \sum_{\chi\in \co}
    	\chi(n)\ \sumstar_{v,w \, (\text{mod}\ q)} \chi\eta_1(w) \eta_1\overline{\eta_2}(v) e_q(vw+v) \\
    	&=  \frac{1}{q|\co|}\ \sumstar_{w \,  (\text{mod}\ q)} \eta_1(w) \sum_{\chi\in \co}
    	\chi(nw)\ \sumstar_{v \, (\text{mod}\ q)} \eta(v) e_q(v(w+1)) \\
    	&=  \frac{1}{q}\ \sumstar_{w \, (\text{mod}\ q)} 
    	\eta_1(w)  \mathbb{E}_\co[\chi(nw)]c_\eta(w+1) .
    \end{align*}
    
    By Lemma \ref{lem:generalized_Gauss_sum}, $c_\eta(w+1) \neq 0$ only if $\ord_p(w+1)=k-h$. By Lemma \ref{lem:full-orbit-char-avg}, $\mathbb{E}_\co[\chi(nw)]\neq 0$ only if $(nw)^{p-1}\equiv 1\ (\text{mod}\ p^{k-1})$, or equivalently $nw\equiv \zeta\ (\text{mod}\ p^{k-1})$ for some $\zeta\in \boldsymbol{\mu}_{p-1}$. Therefore
    $$
    k-h=\ord_p(w+1) = \ord_p(\zeta\bar{n}+1)=\ord_p(\zeta+n).
    $$
    It follows that $n\equiv -\zeta\ (\text{mod}\ p^{k-h})$, and thus  
    $$n^{p-1} \equiv 1 \ (\text{mod}\ p^{k-h}).$$
    The first statement is proved. 
    
    Next we assume that $n^{p-1} \equiv 1 \ (\text{mod}\ p^{k-h})$. We have 
    \begin{align*}
    \frac{1}{|\co|}  \sum_{\chi\in \co} \epsilon_{\chi\eta_1}\epsilon_{\overline{\chi\eta_2}} \chi(n)
    	&= \frac{1}{q}\ \sumstar_{w \, (\text{mod}\ q)}  \eta_1(w)
    	\mathbb{E}_\co[\chi(nw)] c_\eta(w+1)  \\
    	&=  \frac{1}{q} \sum_{\zeta \in \boldsymbol{\mu}_{p-1}} \ 
    	\sumstar_{\substack{w \, (\text{mod}\ q)\\ w\equiv \zeta\bar{n} \, (\text{mod}\ p^{k-1})}} 
    	\eta_1(w) \mathbb{E}_\co[\chi(nw)] c_\eta(w+1) .
    \end{align*}
    By Lemma \ref{lem:generalized_Gauss_sum}, $|c_\eta(w+1)| \leq p^{k-h/2}$. 
    Since $| \mathbb{E}_\co[\chi(nw)]|\leq 1$, $| \eta_1(w)|=1$, and the number of summations is finite ($p$ is a fixed prime), the second statement follows.
\end{proof}

\subsection{Character averages for thin Galois orbits}

In this section, we study character averages for thin Galois orbits.

\begin{lem}\label{lem:thin-orbit-char-avg}
Let $0\leq\kappa\leq k-1$ and $n\in (\Z/q\Z)^\times$. If $\mathbb{E}_{\co_\kappa}[\chi(n)]\neq 0$, then $n^{p-1} \equiv 1 \ (\emph{mod}\ p^{\min(\kappa+1,k-1)})$. 
\end{lem}

\begin{proof}
This is \cite[Lemma 3.1]{KMN16}.
\end{proof}

We have an analog of Lemma \ref{lem:averageing_character_twist_chieta1_barchieta2} for thin Galois orbits.

\begin{lem}\label{lem:averageing_character_twist_chieta1_barchieta2_thin}
    Let $n\in (\Z/q\Z)^\times$. 
    
    \begin{enumerate}
        \item[{\rm (i)}] If $n^{p-1} \not\equiv 1 \ (\emph{mod}\ p^{\min(\kappa+1,k-h)})$, then
    $$
    \frac{1}{|\co_\kappa|}  \sum_{\chi\in \co_\kappa} \epsilon_{\chi\eta_1}\epsilon_{\overline{\chi\eta_2}} \chi(n) = 0.
    $$
        \item[{\rm (ii)}] Suppose that $\kappa > k-h/2$.
        If $n^{p-1} \equiv 1 \ (\emph{mod}\ p^{\min(\kappa+1,k-h)})$, then $n^{p-1} \equiv 1 \ (\emph{mod}\ p^{k-h})$ and
    $$
    \frac{1}{|\co_\kappa|}  \sum_{\chi\in \co_\kappa} \epsilon_{\chi\eta_1}\epsilon_{\overline{\chi\eta_2}} \chi(n) \ll p^{k-h/2-\kappa}.
    $$
    \end{enumerate}
\end{lem}

\begin{proof}
    Opening both Gauss sums as in the proof of Lemma \ref{lem:averageing_character_twist_chieta1_barchieta2}, we get
        \begin{equation*}
    	\frac{1}{|\co_\kappa|}  \sum_{\chi\in \co_\kappa} \epsilon_{\chi\eta_1}\epsilon_{\overline{\chi\eta_2}} \chi(n)
    	=  \frac{1}{q}\ \sumstar_{w \, (\text{mod}\ q)} 
    	\eta_1(w)  \mathbb{E}_{\co_\kappa}[\chi(nw)]c_{\eta_1\overline{\eta_2}}(w+1) .
    \end{equation*}
    By Lemma \ref{lem:generalized_Gauss_sum}, $c_{\eta_1\overline{\eta_2}}(w+1) \neq 0$ only if $\ord_p(w+1)=k-h$. By Lemma \ref{lem:thin-orbit-char-avg}, $\mathbb{E}_{\co_\kappa}[\chi(nw)]\neq 0$ only if $(nw)^{p-1}\equiv 1\ (\text{mod}\ p^{\min(\kappa+1,k-1)})$, or equivalently $nw\equiv \zeta\ (\text{mod}\ p^{\min(\kappa+1,k-1)})$ for some $\zeta\in \boldsymbol{\mu}_{p-1}$. 
    It follows that
    $$
    0\equiv w+1 \equiv  \zeta\bar{n}+1 \equiv (\zeta+n)\bar{n}
    \ (\text{mod}\ p^{\min(\kappa+1,k-h)})
    $$
    Hence $n\equiv -\zeta \ (\text{mod}\ p^{\min(\kappa+1,k-h)})$, and thus
    $$n^{p-1} \equiv 1 \ (\text{mod}\ p^{\min(\kappa+1,k-h)}).$$
    Part {\rm (i)} is proved. 
    
    Next we assume that $\kappa > k-h/2$ and that $n^{p-1} \equiv 1 \ (\text{mod}\ p^{\min(\kappa+1,k-h)})$. It is plain that $\kappa+1 > k-h$, whence $n^{p-1} \equiv 1 \ (\text{mod}\ p^{k-h})$. We have
    \begin{align*}
    \frac{1}{|\co_\kappa|}  \sum_{\chi\in \co_\kappa} \epsilon_{\chi\eta_1}\epsilon_{\overline{\chi\eta_2}} \chi(n)
    &= \frac{1}{q}\ \sumstar_{w \, (\text{mod}\ q)}  \eta_1(w)
    	\mathbb{E}_{\co_\kappa}[\chi(nw)] c_\eta(w+1)  \\
    &= \frac{1}{q} \sum_{\zeta \in \boldsymbol{\mu}_{p-1}} 
    	\sumstar_{\substack{w \, (\text{mod}\ q) \\ 
        w\equiv \zeta\bar{n} \, (\text{mod}\ p^{\kappa+1}) \\
        w+1\equiv 0 \, (\text{mod}\ p^{k-h}) \\
        }} 
    	\eta_1(w) \mathbb{E}_{\co_\kappa}[\chi(nw)] c_\eta(w+1) .
    \end{align*}
    By Lemma \ref{lem:generalized_Gauss_sum}, we have $|c_\eta(w+1)| \leq p^{k-h/2}$. 
    Since  $|\eta_1(w)| = 1$, $|\mathbb{E}_{\co_\kappa}[\chi(nw)]|\leq 1$ and the number of summations is $O(q/p^{\kappa+1})$ ($p$ is a fixed prime), it follows that
    $$
    \frac{1}{|\co_\kappa|}  \sum_{\chi\in \co_\kappa} \epsilon_{\chi\eta_1}\epsilon_{\overline{\chi\eta_2}} \chi(n)
    \ll p^{-h/2} \frac{q}{p^{\kappa+1}} \ll p^{k-h/2-\kappa}.
    $$
    Part {\rm (ii)} is proved. 
\end{proof}

\subsection{Roots of unity and the Roth-Ridout theorem}

In this section, we gather some $p$-adic congruence equations that can be treated elementarily or using more powerful results from $p$-adic Diophantine approximation.

Let $p$ be a prime number. The equation $x^{p-1}-1=0$ has $p-1$ roots in $\Z_p$; these roots form the set $\boldsymbol{\mu}_{p-1}$ of $(p-1)^{\text{th}}$ roots of unity in $\overline{\Q}_p$. For each $1\leq j\leq p-1$, there exists $\zeta_j\in \boldsymbol{\mu}_{p-1}$ such that $\zeta_j\equiv j \ (\text{mod}\ p)$ or equivalently $\zeta_j-j\in p\Z_p$. If $x\in\Z$ satisfies $x^{p-1}\equiv 1 \ (\text{mod}\ p^a)$ for some positive integer $a$, then $x \equiv \zeta_j \ (\text{mod}\ p^a)$, or equivalently $x-\zeta_j\in p^a\Z_p$, for exactly one $1\leq j\leq p-1$.

We recall several results from \cite{KMN16} for later use; Lemmas \ref{lem:Roth-Ridout-a} and \ref{lem:Roth-Ridout-b} below are consequences of the Roth-Ridout theorem.

\begin{lem}\label{lem:congruences-prime-power}
   If $m$ and $\alpha\geq 1$ are integers satisfying $m^{p}\equiv 1 \ (\emph{mod}\ p^\alpha)$, then $m\equiv 1 \ (\emph{mod}\ p^{\alpha-1})$.
\end{lem}

\begin{proof}
    This is \cite[Lemma 2.2]{KMN16}.
\end{proof}

\begin{lem}\label{lem:Roth-Ridout-a}
    Suppose that $0<\delta<\tfrac{1}{2}$ and $\beta\in\Z_p$ is algebraic over $\Q$ of degree at least $2$. Then, for all sufficiently large positive integers $\alpha\geq \alpha_0(\beta,\delta)$, there are no nonzero integers $a$ and $b$ such that $|a|,|b|<(p^\alpha)^{1/2-\delta}$ and that
    $$
    a \equiv \beta b \ (\emph{mod}\ p^{\alpha-1}).
    $$
\end{lem}

\begin{proof}
    This is \cite[Lemma 2.4]{KMN16}.
\end{proof}

\begin{rem}
    It is plain that Lemma \ref{lem:Roth-Ridout-a} also holds if we replace the congruence $a \equiv \beta b \ (\text{mod}\ p^{\alpha-1})$ by the congruence $a \equiv \beta b \ (\text{mod}\ p^{\alpha})$.
\end{rem}

\begin{lem}\label{lem:Roth-Ridout-b}
    Suppose that $0<\delta<\tfrac{1}{2}$ and $\alpha\geq \alpha_1(\delta)$ be sufficiently large. Let $\mathcal{A}_\alpha, \mathcal{B}_\alpha \subset \Z$ be intervals of length at most $(p^\alpha)^{1/2-\delta}$. Then there are at most $(p-1)$ pairs $(a,b)\in \mathcal{A}_\alpha\times \mathcal{B}_\alpha$ such that $(ab,p)=1$, $a \not\equiv \pm b \ (\emph{mod}\ p^{\alpha-1})$ and that
    $$
    a^{p-1} \equiv b^{p-1} \ (\emph{mod}\ p^{\alpha-1}).
    $$
\end{lem}

\begin{proof}
    This is \cite[Lemma 2.5]{KMN16}.
\end{proof}

\begin{rem}
    It is plain that Lemma \ref{lem:Roth-Ridout-b} also holds if we replace the congruence $a^{p-1} \equiv b^{p-1} \ (\text{mod}\ p^{\alpha-1})$ by the congruence $a^{p-1} \equiv b^{p-1} \ (\text{mod}\ p^{\alpha})$.
\end{rem}

For $\alpha\in\mathbb{N}$, define
\begin{align}
D_0(A,B; p^\alpha) &= \frac{1}{\sqrt{AB}} 
\sumstar_{\substack{A\le a < 2A,   B\le b < 2B\\ a^{p-1} \equiv b^{p-1} \ (\text{mod}\ p^{\alpha}) \\ a\neq b}} 1 
= \frac{1}{\sqrt{AB}}\sum_{\zeta \in \boldsymbol{\mu}_{p-1}} \
\sumstar_{\substack{A\le a < 2A,  B\le b < 2B\\ a\equiv \zeta b \ (\text{mod}\ p^{\alpha}) \\ a\neq b}} 1 , \label{eq:dyadic-sum-congruence} \\
D_1(A,B; p^\alpha) &= \frac{1}{\sqrt{AB}} 
\sumstar_{\substack{A\le a < 2A,  B\le b < 2B\\ a\equiv \pm b \ (\text{mod}\ p^{\alpha})\\ a\neq b}} 1 , \label{eq:dyadic-sum-congruence-trivial} \\
D_2(A,B; p^\alpha) &= \frac{1}{\sqrt{AB}}
\sum_{\zeta \in \boldsymbol{\mu}_{p-1}\setminus \{\pm 1\}} \
\sumstar_{\substack{A\le a < 2A, B\le b < 2B \\ a\equiv \zeta b \ (\text{mod}\ p^{\alpha})}} 1  \label{eq:dyadic-sum-congruence-RR} .
\end{align}
These are the sums that we need to estimate in the moment computations. Note that $$D_0(A,B; p^\alpha)=D_1(A,B; p^\alpha)+D_2(A,B; p^\alpha).$$
The sums $D_0(A,B; p^\alpha)$ and $D_1(A,B; p^\alpha)$ can be treated elementarily, whereas $D_2(A,B; p^\alpha)$ can be estimated by the Roth-Ridout theorem.

\begin{lem}\label{lem:dyadic-congruence-naive}
    Let $\alpha\in\mathbb{N}$. 
We have
$$
D_0(A,B; p^\alpha) \ll \frac{\sqrt{AB}}{p^\alpha} + \min\left(\sqrt{\frac{A}{B}},\sqrt{\frac{B}{A}}\right) .
$$
In particular, $D_0(A,B; p^\alpha) = O(\sqrt{AB}/p^\alpha + 1)$.
\end{lem}

\begin{proof}
    Without loss of generality, assume that $B\leq A$. For each $B\leq b<2 B$, there are $O(A/p^\alpha+1)$ values of $A\leq a<2A$ with $a\equiv \zeta b\ (\text{mod}\ p^{\alpha})$ for each $\zeta \in \boldsymbol{\mu}_{p-1}$. Then 
$$
D_0(A,B; p^\alpha) \ll \frac{1}{\sqrt{AB}}\Big(\frac{AB}{p^\alpha}+B\Big)
= \frac{\sqrt{AB}}{p^\alpha} + \sqrt{\frac{B}{A}}.$$
The lemma is proved.
\end{proof}

\begin{lem}\label{lem:dyadic-congruence-naive-1}
    Let $\alpha\in\mathbb{N}$. 
We have
$$
D_1(A,B; p^\alpha) \ll \frac{\sqrt{AB}}{p^\alpha}.
$$
\end{lem}

\begin{proof}
    If $A\le a < 2A$ and $B\le b < 2B$ satisfy $a\equiv \pm b \ (\text{mod}\ p^{\alpha})$ and $a\neq b$, then either $a\gg p^{\alpha}$ or $b\gg p^{\alpha}$. Without loss of generality, assume that $a\gg p^{\alpha}$, whence $\sqrt{B/A} \ll \sqrt{AB}/p^\alpha$. By Lemma \ref{lem:dyadic-congruence-naive}, we have 
    $$D_1(A,B; p^\alpha) \leq D_0(A,B; p^\alpha)
    \ll \frac{\sqrt{AB}}{p^\alpha}  + \sqrt{\frac{B}{A}}
    \ll \frac{\sqrt{AB}}{p^\alpha},$$
as claimed.
\end{proof}

We are in a position to derive the main estimates for the number of solutions of the $p$-adic congruence equations that naturally arise from Dirichlet characters in Galois orbit. Although its proof uses arguments that appeared in \cite{KMN16}, Proposition \ref{prop:dyadic-congruence} as formulated below is new.

\begin{prop}\label{prop:dyadic-congruence}
    Suppose that $0<\delta<\tfrac{1}{4}$ and that $\alpha\geq \alpha_2(\delta)$ is a sufficiently large positive integer. 
    \begin{enumerate}
        \item[{\rm (i)}] If $2A < p^{\alpha(1/2-\delta)}$ and $2B < p^{\alpha(1/2-\delta)}$, then 
        $$
        D_0(A,B;p^\alpha)=D_1(A,B;p^\alpha)=D_2(A,B;p^\alpha)=0.
        $$
        \item[{\rm (ii)}] If $2A \geq p^{\alpha(1/2-\delta)}$ and $2B \geq p^{\alpha(1/2-\delta)}$, then 
        $$
        D_0(A,B;p^\alpha) \ll \frac{\sqrt{AB}}{p^{\alpha}} \cdot p^{4\delta\alpha} . 
        $$
        \item[{\rm (iii)}] If $2A \geq p^{\alpha(1/2-\delta)}$ and $2B < p^{\alpha(1/2-\delta)}$, then $D_2(A,B;p^\alpha) \ll \frac{\sqrt{A/B}}{p^{\alpha/2}}$ and
        $$
        D_0(A,B;p^\alpha) \ll \frac{\sqrt{AB}}{p^{\alpha}} + 
        p^{\alpha(-1/4+\delta)} . 
        $$
    \end{enumerate}
    In particular, in all cases we have 
    $$
        D_0(A,B;p^\alpha) \ll \frac{\sqrt{AB}}{p^{\alpha}} \cdot p^{4\delta\alpha} + 
        p^{\alpha(-1/4+\delta)} . 
    $$
\end{prop}

\begin{proof}
    Statement ${\rm (i)}$ is an immediate consequence of Lemma \ref{lem:Roth-Ridout-a}.

    To prove ${\rm (ii)}$, let $\zeta \in \boldsymbol{\mu}_{p-1}\setminus \{\pm 1\}$ and divide the dyadic intervals $A\leq a<2A$ and $B\leq b < 2B$ into sub-intervals of the same length $Q:=p^{\alpha(1/2-2\delta)}$. 
    For each pair of those sub-intervals, the number of pairs $(a,b)$ in this box and satisfying $a\equiv \zeta b\ (\text{mod}\ p^{\alpha})$ is $O(1)$, by the Roth-Ridout theorem (Lemma \ref{lem:Roth-Ridout-b}). The contribution of those terms to the sum gives
    $$
    D_2(A,B;p^\alpha) \ll \frac{1}{\sqrt{AB}} \cdot \frac{AB}{Q^2} \ll \frac{\sqrt{AB}}{p^{\alpha(1-4\delta)}}.
    $$
    This estimate together with Lemma \ref{lem:dyadic-congruence-naive-1} yield ${\rm (ii)}$.
    
    The proof of ${\rm (iii)}$ is similar. Let $\zeta \in \boldsymbol{\mu}_{p-1}\setminus \{\pm 1\}$ and divide the dyadic interval $A\leq a < 2A$ into sub-intervals of the same length $Q$ as above. For each of those sub-intervals, the number of pairs $(a,b)$ with $a$ in the sub-interval and $B\leq b < 2B$ satisfying $a\equiv \zeta b\ (\text{mod}\ p^{\alpha})$ is $O(1)$, by the Roth-Ridout theorem (Lemma \ref{lem:Roth-Ridout-b}). The contribution of those terms to the sum gives
    $$
    D_2(A,B;p^\alpha) \ll \frac{1}{\sqrt{AB}}  \cdot  \frac{A}{Q} \ll \frac{\sqrt{A/B}}{p^{\alpha(1/2-2\delta)}}.
    $$
    This estimate together with Lemma \ref{lem:dyadic-congruence-naive-1} yield 
    \begin{equation}\label{eq:dyadic-congruence-balanced1}
        D_0(A,B;p^\alpha) \ll \frac{\sqrt{AB}}{p^{\alpha}} + \frac{\sqrt{A/B}}{p^{\alpha/2}} \cdot p^{2\delta\alpha} . 
    \end{equation}

    Let $\lambda>0$ be a parameter to be chosen.
    For the `unbalanced' range 
    $$
    \frac{A}{B}\notin (p^{-\lambda},p^{\lambda}),
    $$
    by assumption we have $\frac{A}{B}\geq p^\lambda$. It follows from Lemma \ref{lem:dyadic-congruence-naive} that 
    \begin{equation}\label{eq:dyadic-congruence-unbalanced}
    D_0(A,B;p^\alpha) \ll \frac{\sqrt{AB}}{p^\alpha} + p^{-\lambda/2}. 
    \end{equation}
    In the `balanced' range 
    $$ p^{-\lambda} < \frac{A}{B} < p^{\lambda},$$ 
    it follows from \eqref{eq:dyadic-congruence-balanced1} that
    \begin{equation}\label{eq:dyadic-congruence-balanced2}
        D_0(A,B;p^\alpha) \ll \frac{\sqrt{AB}}{p^{\alpha}} + p^{(\lambda-\alpha)/2} \cdot p^{2\delta\alpha} . 
    \end{equation}
    On combining \eqref{eq:dyadic-congruence-unbalanced} and \eqref{eq:dyadic-congruence-balanced2} and setting $p^\lambda=p^{\alpha(1/2-2\delta)}$, we conclude ${\rm (iii)}$.
\end{proof}


\section{Twisted moments}

We set to compute the twisted moments over full Galois orbits and thin Galois orbits.

\subsection{Twisted moments over full Galois orbits}

In this section, we prove Theorem \ref{thm:twisted-moment-full-orbit} together with its consequences (Corollaries \ref{cor:full-orbit-moment} and \ref{cor:full-orbit-nonvanishing}).

\begin{proof}[Proof of Theorem \ref{thm:twisted-moment-full-orbit}]
Let
$$
S(m_1,m_2) 
    = \frac{1}{|\co|}\sum_{\substack{\chi\in\co}}L(\tfrac12,\chi\eta_1)L(\tfrac12,\overline{\chi\eta_2})\chi(m_1)\overline{\chi}(m_2).
$$
By Lemma \ref{lem:afe}, we have 
\begin{equation*}
\label{eq:funceq2ndmol}
S(m_1,m_2)= S^+(m_1,m_2)+ i^{-(\kappa_1+\kappa_2)}S^-(m_1,m_2),   
\end{equation*}
where
\begin{equation}
    \label{eq:funceq2ndmolplus}
S^+(m_1,m_2) =
    \sum_{\substack{n_1,n_2\ge 1}} 
    \frac{\eta_1(n_1)\overline{\eta_2}(n_2)}{\sqrt{n_1n_2}}V\Big(\frac{\pi n_1n_2}{qX}\Big) 
    \mathbb{E}_\co[\chi(m_1n_1\overline{m_2n_2})]  
\end{equation}
and 
\begin{equation}
\label{eq:funceq2ndmolminus}
 S^-(m_1,m_2)= 
     \sum_{\substack{n_1,n_2\ge 1}} 
    \frac{\overline{\eta_1}(n_1)\eta_2(n_2)}{\sqrt{n_1n_2}}V\Big(\frac{\pi n_1n_2X}{q}\Big) \frac{1}{|\co|} 
    \sum_{\chi\in \co} \epsilon_{\chi\eta_1} \epsilon_{\overline{\chi}\overline{\eta_2}} 
    \chi(m_1\overline{n_1}\overline{m_2}n_2).
\end{equation}

For $S^+(m_1,m_2)$ we write 
\[
S^+(m_1,m_2)=M(m_1,m_2)+E^+(m_1,m_2),
\]
where $M(m_1,m_2)$ is the contribution of the diagonal terms $m_1n_1=m_2n_2$, and $E^+(m_1,m_2)$ is the contribution of the off-diagonal terms $m_1n_1\ne m_2n_2$.  Write $m_1=gm_1',m_2=gm_2'$, where $(m_1',m_2')=1$. Then
\begin{align}\label{Mformula}
M(m_1,m_2)=\sum_{\substack{n_1,n_2\geq 1\\m_1'n_1=m_2'n_2}}\frac{\eta_1(n_1)\overline{\eta_2}(n_2)}{\sqrt{n_1n_2}}V\Big(\frac{\pi n_1n_2}{qX}\Big)=\frac{\eta_1(m_2')\overline{\eta_2}(m_1')}{\sqrt{m_1'm_2'}}\sum_{n\geq 1}\frac{\eta(n)}{n}V\Big(\frac{\pi m_1'm_2'n^2}{qX}\Big).
\end{align}
By \eqref{formulaV+}, we get that
\begin{align*}
M(m_1,m_2) =&\frac{\eta_1(m_2')\overline{\eta_2}(m_1')}{\sqrt{m_1'm_2'}}\frac{1}{2\pi i}\int_{(1)}\frac{\Gamma\big(\frac{1/2+s+\kappa_1}{2}\big)\Gamma\big(\frac{1/2+s+\kappa_2}{2}\big)}{\Gamma\big(\frac{1/2+\kappa_1}{2}\big)\Gamma\big(\frac{1/2+\kappa_2}{2}\big)}G(s)\Big(\frac{qX}{\pi m_1'm_2'}\Big)^{s}L(1+2s,\eta)\frac{ds}{s}.
\end{align*}
We move the contour to Re$(s)=-1/2+\varepsilon$, crossing a simple pole at $s=0$. Estimating the new integral trivially gives 
\begin{equation}\label{eq:diagonalterms}
M(m_1,m_2) =\frac{\eta_1(m_2')\overline{\eta_2}(m_1')}{\sqrt{m_1'm_2'}}L(1,\eta)+O_\varepsilon(p^{h/2}q^{-1/2+\varepsilon}X^{-1/2}).
\end{equation}

\kommentar{For the application to Theorem \ref{thm:nonvanishing-full-orbit}, it is more convenient to have the main term in the form of a truncated Euler product. To this end, we use \eqref{formulaV+} to write \eqref{Mformula} as
\begin{align*}
&\frac{\eta_1(m_2')\overline{\eta_2}(m_1')}{\sqrt{m_1'm_2'}}\frac{1}{2\pi i}\int_{(1)}\frac{\Gamma\big(\frac{1/2+s+\kappa_1}{2}\big)\Gamma\big(\frac{1/2+s+\kappa_2}{2}\big)}{\Gamma\big(\frac{1/2+\kappa_1}{2}\big)\Gamma\big(\frac{1/2+\kappa_2}{2}\big)}G(s)\Big(\frac{qX}{\pi m_1'm_2'}\Big)^{s}\sum_{n\geq 1}\frac{\eta(n)}{n^{1+2s}}\frac{ds}{s}.
\end{align*}
We can truncate the sum to $n\leq q^{2}$ since the contribution of the terms $n>q^2$ can be seen to be negligible by moving the contour far to the right. By the same argument, we can then extend the sum to all $q^2$-smooth $n$. Thus, with a negligible error term, 
\begin{align}\label{newM}
M(m_1,m_2)&=\frac{\eta_1(m_2')\overline{\eta_2}(m_1')}{\sqrt{m_1'm_2'}}\\
&\qquad\times\frac{1}{2\pi i}\int_{(1)}\frac{\Gamma\big(\frac{1/2+s+\kappa_1}{2}\big)\Gamma\big(\frac{1/2+s+\kappa_2}{2}\big)}{\Gamma\big(\frac{1/2+\kappa_1}{2}\big)\Gamma\big(\frac{1/2+\kappa_2}{2}\big)}G(s)\Big(\frac{qX}{\pi m_1'm_2'}\Big)^{s}\prod_{p\leq q^2}\bigg(1-\frac{\eta(p)}{p^{1+2s}}\bigg)^{-1}\frac{ds}{s}.\nonumber
\end{align}}
To estimate 
$$
\mathcal{E}^+(\vartheta):=\sum_{m_1,m_2\leq q^\vartheta}\frac{x_{m_1}\overline{x_{m_2}}}{\sqrt{m_1m_2}}E^+(m_1,m_2),    
$$
we reason as in Section $3.3.2$ of \cite{KMN16}.
First, it follows from \eqref{eq:funceq2ndmolplus} and Lemma \ref{lem:full-orbit-char-avg} that
$$
\mathcal{E}^+(\vartheta) \ll_\varepsilon  q^\varepsilon 
\sumstar_{\substack{m_1,m_2 \leq q^\vartheta\\n_1n_2\leq q^{1+\varepsilon}X\\ (m_1n_1)^{p-1} \equiv (m_2n_2)^{p-1} \ (\text{mod}\ p^{k-1})\\ m_1n_1\neq m_2n_2}} \frac{1}{\sqrt{m_1m_2n_1n_2}}.
$$
Put $a=m_1n_1$ and $b=m_2n_2$ and decompose the sum into dyadic intervals. This reduces to $O_\varepsilon(q^{\varepsilon})$ sums of the form $D_0(A,B;p^{k-1})$ with $1\leq  AB\leq q^{1+2\vartheta+\varepsilon}X$. Therefore
\begin{equation}
\label{eq:dyadicS+} 
\mathcal{E}^+(\vartheta) \ll_\varepsilon  q^\varepsilon 
\max_{\substack{A\geq 1,B\geq 1 \\  AB\leq q^{1+2\vartheta+\varepsilon}X }}
D_0(A,B;p^{k-1}) .
\end{equation}
We apply Proposition \ref{prop:dyadic-congruence}, setting $\delta=\varepsilon$ therein, to derive 
\begin{equation}\label{eq:S+errorterm-condition}
\mathcal{E}^+(\vartheta)  \ll_\varepsilon
q^{-1/2+\vartheta+\varepsilon}X^{1/2} 
+ q^{-1/4+\varepsilon}.
\end{equation}

We next move on to $S^-(m_1,m_2)$. Contrary to $S^+(m_1,m_2)$, in $S^-(m_1,m_2)$ there is no main term but only the error term
$$
\mathcal{E}^-(\vartheta):=\sum_{m_1,m_2\leq q^\vartheta}\frac{x_{m_1}\overline{x_{m_2}}}{\sqrt{m_1m_2}}S^-(m_1,m_2).
$$ 
From \eqref{eq:funceq2ndmolminus} and Lemma \ref{lem:averageing_character_twist_chieta1_barchieta2}, we deduce that
\begin{align*}
\mathcal{E}^-(\vartheta) \ll_\varepsilon q^\varepsilon p^{-h/2}\ 
\sumstar_{\substack{m_1,m_2 \leq q^\vartheta\\n_1n_2\leq q^{1+\varepsilon}/X\\ (m_1n_2)^{p-1} \equiv (m_2n_1)^{p-1} \, (\text{mod}\ p^{k-h}) \\ m_1n_2\neq m_2n_1}} \frac{1}{\sqrt{m_1m_2n_1n_2}}.
\end{align*}
Put $a=m_1n_2$ and $b=m_2n_1$ and decompose the sum into dyadic intervals. 
This reduces to $O_\varepsilon(q^{\varepsilon})$ sums of the form $D_0(A,B;p^{k-h})$
with $1\leq AB\leq q^{1+2\vartheta+\varepsilon}/X$. Therefore
\begin{equation*}
\mathcal{E}^-(\vartheta)  \ll_\varepsilon q^\varepsilon p^{-h/2}
\max_{\substack{A\geq 1,B\geq 1 \\  AB\leq q^{1+2\vartheta+\varepsilon} / X }}
D_0(A,B;p^{k-h}) .
\end{equation*}
It follows from Proposition \ref{prop:dyadic-congruence} (with the choice $\delta=\varepsilon$) that
\begin{align}
\mathcal{E}^-(\vartheta)&  \ll_\varepsilon q^\varepsilon p^{-h/2}
\max_{\substack{A\geq 1,B\geq 1 \\  AB\leq q^{1+2\vartheta+\varepsilon} / X }} \frac{\sqrt{AB}}{p^{k-h}}  + p^{-h/4}q^{-1/4+\varepsilon} \nonumber \\
& \ll_\varepsilon 
p^{h/2}q^{-1/2+\vartheta+\varepsilon}X^{-1/2} + p^{-h/4}q^{-1/4+\varepsilon}. 
\end{align}

For all $0<h<k$, it follows from Lemma \ref{lem:dyadic-congruence-naive} that

\begin{align}
\mathcal{E}^-(\vartheta)&  \ll_\varepsilon q^\varepsilon p^{-h/2}
\max_{\substack{A\geq 1,B\geq 1 \\  AB\leq q^{1+2\vartheta+\varepsilon} / X }} \frac{\sqrt{AB}}{p^{k-h}} +  q^\varepsilon p^{-h/2} \nonumber \\
& \ll_\varepsilon 
p^{h/2}q^{-1/2+\vartheta+\varepsilon}X^{-1/2} +  p^{-h/2}q^\varepsilon. \label{eq:boundS-zero}
\end{align}
Assuming further that $h<(1-\varepsilon)k$, we can apply Proposition \ref{prop:dyadic-congruence} (with the choice $\delta=\varepsilon$) to obtain
\begin{equation}
\mathcal{E}^-(\vartheta)  \ll_\varepsilon p^{h/2}q^{-1/2+\vartheta+\varepsilon}X^{-1/2} + p^{-h/4}q^{-1/4+\varepsilon}.
\label{eq:boundS-}
\end{equation}

In view of \eqref{eq:diagonalterms}, \eqref{eq:S+errorterm-condition}, \eqref{eq:boundS-zero} and \eqref{eq:boundS-} we derive that
$$
S(m_1,m_2)=\frac{\eta_1(m_2')\overline{\eta_2}(m_1')}{\sqrt{m_1'm_2'}}L(1,\eta)+E(m_1,m_2),
$$
where
\begin{align*}
\mathcal{E}(\vartheta):=\sum_{m_1,m_2\leq q^\vartheta}\frac{x_{m_1}\overline{x_{m_2}}}{\sqrt{m_1m_2}}E(m_1,m_2)
= \mathcal{E}^+(\vartheta) + \mathcal{E}^-(\vartheta)
\end{align*}
satisfies the following property. We have that for all $0<h<k$, 
\begin{align*}
\mathcal{E}(\vartheta) \ll_\varepsilon q^{-1/2+\vartheta+\varepsilon}X^{1/2}
+ p^{h/2} q^{-1/2+\vartheta+\varepsilon}X^{-1/2}
+ p^{-h/2}q^{\varepsilon} + q^{-1/4+\varepsilon}
\end{align*}
and for $0<h<(1-\varepsilon)k$,
\begin{align*}
\mathcal{E}(\vartheta)\ll_\varepsilon q^{-1/2+\vartheta+\varepsilon}X^{1/2}
+ p^{h/2} q^{-1/2+\vartheta+\varepsilon}X^{-1/2}
+ p^{-h/4}q^{-1/4+\varepsilon}.
\end{align*}
Theorem \ref{thm:twisted-moment-full-orbit} follows by choosing $X=p^{h/2}$.

\end{proof}


We are in a position to derive two consequences of Theorem \ref{thm:twisted-moment-full-orbit}.

\begin{proof}[Proof of Corollary \ref{cor:full-orbit-moment}]
    This corollary follows immediately from Theorem \ref{thm:twisted-moment-full-orbit} on setting $m_1=m_2=1$.
\end{proof}

\begin{proof}[Proof of Corollary \ref{cor:full-orbit-nonvanishing}]
   This is a consequence of Corollary \ref{cor:full-orbit-moment} and the subconvexity bound for Dirichlet $L$-functions. From Corollary \ref{cor:full-orbit-moment} and  Siegel's Theorem (see, for instance, \cite[Theorem 11.14]{MV07}), it follows that
\begin{align}\label{cor:full-orbit-nonvanishing for1}
&\sum_{\substack{\chi\in\co}}L(\tfrac12,\chi\eta_1)L(\tfrac12,\overline{\chi\eta_2})\gg_\varepsilon q^{1-\varepsilon}.
\end{align}
By the Weyl bound $L(1/2,\chi\eta_1)\ll_\varepsilon q^{1/6+\varepsilon}$ \cite[Theorem 1.1]{PY23}, we also have that the left hand side is 
\begin{align}\label{cor:full-orbit-nonvanishing for2}
\ll_\varepsilon q^{1/3+\varepsilon}\cdot \#\big\{\chi\in\co:L(\tfrac12,\chi\eta_1)L(\tfrac12,\chi\eta_2)\ne0\big\}.
\end{align}
Combining \eqref{cor:full-orbit-nonvanishing for1} and \eqref{cor:full-orbit-nonvanishing for2} we obtain Corollary \ref{cor:full-orbit-nonvanishing}.
\end{proof}

\subsection{Twisted moments over thin Galois orbits}

In this section, we prove Theorem \ref{thm:twisted-moment-thin-orbit}.

\begin{proof}[Proof of Theorem \ref{thm:twisted-moment-thin-orbit}]

We proceed as in the proof of Theorem \ref{thm:twisted-moment-full-orbit}, with Lemmas \ref{lem:thin-orbit-char-avg} and \ref{lem:averageing_character_twist_chieta1_barchieta2_thin} in place of 
Lemmas \ref{lem:full-orbit-char-avg} and \ref{lem:averageing_character_twist_chieta1_barchieta2} respectively.

    Let
$$
S_\kappa(m_1,m_2) 
    = \frac{1}{|\co_\kappa|}\sum_{\substack{\chi\in\co_\kappa}}L(\tfrac12,\chi\eta_1)L(\tfrac12,\overline{\chi\eta_2})\chi(m_1)\overline{\chi}(m_2).
$$
By Lemma \ref{lem:afe}, we have 
\begin{equation*}
\label{eq:funceq2ndmol-thin}
S_\kappa(m_1,m_2)= S_\kappa^+(m_1,m_2)+ i^{-(\kappa_1+\kappa_2)}S_\kappa^-(m_1,m_2),   
\end{equation*}
where
\begin{equation}
\label{eq:funceq2ndmolplus-thin}
S_\kappa^+(m_1,m_2) =
    \sum_{\substack{n_1,n_2\ge 1}} 
    \frac{\eta_1(n_1)\overline{\eta_2}(n_2)}{\sqrt{n_1n_2}}V\Big(\frac{\pi n_1n_2}{qX}\Big) \mathbb{E}_{\co_\kappa}[\chi(m_1n_1\overline{m_2n_2})]   
\end{equation}
and 
\begin{equation}
\label{eq:funceq2ndmolminus-thin}
 S_\kappa^-(m_1,m_2)= 
     \sum_{\substack{n_1,n_2\ge 1}} 
    \frac{\eta_1(\overline{n_1})\eta_2(n_2)}{\sqrt{n_1n_2}}V\Big(\frac{\pi n_1n_2X}{q}\Big) \frac{1}{|\co_\kappa|} 
    \sum_{\chi\in \co_\kappa} \epsilon_{\chi\eta_1} \epsilon_{\overline{\chi}\overline{\eta_2}} 
    \chi(m_1\overline{n_1}\overline{m_2}n_2).
\end{equation}

For $S_\kappa^+(m_1,m_2)$ we write  
\[
S_\kappa^+(m_1,m_2)=M_\kappa(m_1,m_2)+E_\kappa^+(m_1,m_2),
\]
where $M_\kappa(m_1,m_2)$ is the contribution of the diagonal terms $m_1n_1=m_2n_2$, and $E_\kappa^+(m_1,m_2)$ is the contribution of the off-diagonal terms $m_1n_1\ne m_2n_2$.  
The main term $M_\kappa(m_1,m_2)$ can be evaluated asymptotically similarly to the computation of $M(m_1,m_2)$. 
We note that 
$$M_\kappa(m_1,m_2) = M(m_1,m_2),$$ and $M(m_1,m_2)$ is evaluated in \eqref{eq:diagonalterms}. 
We now consider the error term sum
$$
\mathcal{E}^+_\kappa(\vartheta):= 
\sum_{m_1,m_2\leq q^\vartheta}\frac{x_{m_1}\overline{x_{m_2}}}{\sqrt{m_1m_2}}E_\kappa^+(m_1,m_2)
$$
For convenience, put
$$\tilde{\kappa} = \min(\kappa+1,k-1)= \kappa+1.$$
It follows from \eqref{eq:funceq2ndmolplus-thin} and Lemma \ref{lem:thin-orbit-char-avg} that 
$$
\mathcal{E}^+_\kappa(\vartheta)\ll_\varepsilon q^\varepsilon
\max_{\substack{A,B\geq 1\\ 1\leq AB\leq q^{1+2\vartheta+\varepsilon}X}} \
D_0(A,B;p^{\tilde{\kappa}}).
$$
By virtue of Proposition \ref{prop:dyadic-congruence} (with the choice $\delta=\varepsilon$), we deduce that
\begin{equation}\label{eq:S+errorterm-condition-thin}
    \mathcal{E}^+_\kappa(\vartheta)\ll_\varepsilon q^{1/2+\vartheta+\varepsilon}X^{1/2}p^{-\tilde{\kappa}} + p^{-\tilde{\kappa}/4}q^\varepsilon.
\end{equation}

For $S_\kappa^-(m_1,m_2)$, as before there is no main term but only the error term
$$
\mathcal{E}^-_\kappa(\vartheta) := \sum_{m_1,m_2\leq q^\vartheta}\frac{x_{m_1}\overline{x_{m_2}}}{\sqrt{m_1m_2}}S_\kappa^-(m_1,m_2).
$$ 
By \eqref{eq:funceq2ndmolminus-thin}, Lemma \ref{lem:averageing_character_twist_chieta1_barchieta2_thin} (with the assumption $\kappa > k-h/2$) and dyadic decomposition, 
\begin{equation}
\label{eq:dyadicS-thin} 
\mathcal{E}^-_\kappa(\vartheta) \ll_\varepsilon p^{k-\kappa-h/2}q^\varepsilon
\max_{\substack{A,B\geq 1\\ AB\leq q^{1+2\vartheta+\varepsilon}/X}} \
D_0(A,B;p^{k-h}).
\end{equation}

For all $0<h<k$, it follows from Lemma \ref{lem:dyadic-congruence-naive} that 
\begin{equation}\label{eq:boundS-zero-thin}
\mathcal{E}^-_\kappa(\vartheta) \ll_\varepsilon 
p^{h/2-\kappa}q^{1/2+\vartheta+\varepsilon}X^{-1/2} + p^{k-\kappa-h/2}q^\varepsilon.    
\end{equation}
Assuming further that $0<h<(1-\varepsilon)k$, we can apply Proposition \ref{prop:dyadic-congruence}, setting $\delta=\varepsilon$ therein, to obtain
\begin{equation}
\mathcal{E}^-_\kappa(\vartheta) \ll_\varepsilon p^{h/2-\kappa} q^{1/2+\vartheta+\varepsilon}X^{-1/2} + p^{-\kappa+(3k-h)/4}q^{\varepsilon}.
\label{eq:boundS-thin}
\end{equation}

In view of \eqref{eq:S+errorterm-condition-thin}, \eqref{eq:boundS-zero-thin} and \eqref{eq:boundS-thin} we derive that
\[
S_\kappa(m_1,m_2)=\frac{\eta_1(m_2')\overline{\eta_2}(m_1')}{\sqrt{m_1'm_2'}}L(1,\eta)+E_\kappa(m_1,m_2),
\]
where
$$
\mathcal{E}_\kappa(\vartheta) := \sum_{m_1,m_2\leq q^\vartheta}\frac{x_{m_1}\overline{x_{m_2}}}{\sqrt{m_1m_2}}E_\kappa(m_1,m_2) 
= \mathcal{E}^+_\kappa(\vartheta) + \mathcal{E}^-_\kappa(\vartheta)
$$
satisfies the following property. We have that for all $0<h<k$,
$$
\mathcal{E}_\kappa(\vartheta) \ll_\varepsilon 
q^{1/2+\vartheta+\varepsilon}X^{1/2}p^{-\kappa} + p^{-\kappa/4}q^\varepsilon
+ p^{h/2-\kappa} q^{1/2+\vartheta+\varepsilon}X^{-1/2} + p^{k-\kappa-h/2}q^{\varepsilon}
$$
and for $0<h<(1-\varepsilon)k$,
$$
\mathcal{E}_\kappa(\vartheta) \ll_\varepsilon 
q^{1/2+\vartheta+\varepsilon}X^{1/2}p^{-\kappa} + p^{-\kappa/4}q^\varepsilon
+ p^{h/2-\kappa} q^{1/2+\vartheta+\varepsilon}X^{-1/2} + p^{-\kappa+(3k-h)/4}q^{\varepsilon}.
$$
Theorem \ref{thm:twisted-moment-thin-orbit} follows by choosing $X=p^{h/2}$, noting that the latter bound also holds when $(1-\varepsilon)k\leq h < k$.
\end{proof}



\begin{proof}[Proof of Corollary \ref{cor:thin-orbit-moment}]
    This corollary follows immediately from Theorem \ref{thm:twisted-moment-thin-orbit} on setting $m_1=m_2=1$.
\end{proof}


\section{Proof of Theorem \ref{mol_ub}}
\label{mol_section}

Using Theorem $2.1$ in \cite{Cha09}, the following lemma immediately follows.
\begin{lem}
\label{lem_log}
Let $\chi$ be a primitive character modulo $q$. Let $\lambda$ be the solution to $e^{-\lambda} = \lambda+\lambda^2/2$. For any $x$ such that $\log x \geq 2$, we have that
\begin{align*}
\log |L(1/2,\chi)| \leq \Re \sum_{n \leq x} \frac{ \Lambda(n) \chi(n) }{n^{\frac{1}{2}+\frac{\lambda}{\log x}} \log n } \frac{ \log(x/n)}{\log x} + \frac{1+\lambda}{2}\cdot  \frac{ \log q}{\log x} + O \Big( \frac{1}{ \log^2x} \Big).
\end{align*} 
\end{lem}

The Euler product mollifier is defined as follows.

We consider the intervals
$$I_0= (1, q^{\beta_0}], \, I_1=(q^{\beta_0}, q^{\beta_1}], \ldots, I_K= (q^{\beta_{K-1}}, q^{\beta_K}],$$ and for $u \leq K$, let
$$P_{I_j}(\chi; u) = \sum_{p \in I_j} \frac{\chi(p) a(p;u)}{\sqrt{p}} ,$$ where
$$a(p;u) = \Big(1 - \frac{\log p}{\beta_u \log q} \Big) \frac{1}{p^{\frac{\lambda}{\beta_u \log q}}}.  $$
We extend $a(p;u)$ to a completely multiplicative function in the first variable. 

For $t \in \mathbb{R}$ and any integer $\ell$, we define
$$E_{\ell}(t) = \sum_{s \leq \ell} \frac{t^s}{s!}.$$
We note that when $\ell$ is even, then $E_{\ell}(t) >0$. We also have that for $t \leq \ell/e^2$, we have
\begin{equation}
e^t \leq (1+e^{-\ell/2})E_{\ell}(t).
\label{ineq_exp}
\end{equation}
We let $\nu(n)$ be the multiplicative function which is defined by $\nu(p^a) = \frac{1}{a!}$ for prime powers, and $\nu_j(n)$ denotes the $j$-fold convolution of $\nu$; namely, $\nu_j(n) = (\nu * \ldots * \nu)(n)$. Note that $\nu_j(p^a) = \frac{j^a}{a!}$.

We define
\begin{equation}
\label{mj}
M_j(\chi ) := E_{\ell_j}(-P_{I_j}(\chi;K)) = \sum_{\substack{p|n \Rightarrow p \in I_j \\ \Omega(n) \leq \ell_j}} \frac{ \chi(n) \lambda(n) a(n;K) \nu(n)}{\sqrt{n}},
\end{equation} and 
\begin{equation}
\label{mollifier}
M(\chi ) = \prod_{j=0}^K M_j(\chi),
\end{equation}
where $\lambda(n)$ denotes the Liouville function.
We have that 
$$M_j(\chi )^v= \sum_{\substack{p |n \Rightarrow p \in I_j \\ \Omega(n) \leq v \ell_j}}\frac{ \chi (n) \lambda(n) a(n;K) \nu_v(n;\ell_j)}{\sqrt{n}} ,$$
where 
$$ \nu_v(n; \ell_j) = \sum_{\substack{ n = n_1 \cdot \ldots \cdot n_v \\ \Omega(n_i) \leq \ell_j}} \nu(n_1) \cdot \ldots \cdot \nu(n_v).$$

Let 
$$D_{j,v} (\chi ) = \prod_{r=0}^j  (1+e^{-\ell_r/2}) E_{\ell_r}(v \Re P_{I_r}(\chi ;j)),$$ and 
\begin{equation}
S_{j,v}(\chi ) = \exp \Big( v \Re \sum_{p \leq q^{\beta_j/2}} \frac{ \chi (p^2)b(p;j)}{p} \Big),
\label{sjk}
\end{equation} where
\begin{equation}
b(p;j) = \Big(1 - \frac{2 \log p}{\beta_j \log q} \Big) \frac{1}{p^{\frac{2\lambda}{\beta_j \log q}}}.
\label{bp}
\end{equation}
Note that we have 
$$E_{\ell_r} (v \Re P_{I_r} (\chi ;j)) = \sum_{\Omega(mn) \leq \ell_r} \frac{ (v/2)^{\Omega(mn)} \chi(n)  \overline{\chi} (m) a(mn;j) \nu(n) \nu(m)}{\sqrt{mn}} .$$
For $0 \leq r \leq K$, let 
$$\mathcal{T}_r = \Big\{  \chi \pmod q \text{ primitive}, \max_{r \leq u \leq K} |\Re P_{I_r}(\chi ;u) | \leq \frac{\ell_r}{ve^2} \Big\}.$$
Then we have the following proposition.
\begin{prop}
\label{main_prop}
We either have
$$ \max_{0 \leq u \leq K} \Big| \Re P_{I_0}(\chi ;u)\Big|> \frac{\ell_0}{ve^2},$$ or
\begin{align*}
|L(1/2,\chi )|^v  & \ll \exp \Big(  \frac{v (1+\lambda)}{2 \beta_K } \Big) D_{K,v}(\chi ) S_{K,v}(\chi ) \\
&+ \sum_{\substack{0 \leq j \leq K-1 \\ j < u \leq K}} \exp \Big( \frac{ v(1+\lambda) }{2\beta_j } \Big) D_{j,v} (\chi)S_{j,v}(\chi ) \Big(\frac{ve^2 \Re P_{I_{j+1}}(\chi ;u)}{\ell_{j+1}}  \Big)^{s_{j+1}},
\end{align*}
for any even integers $s_{j+1}$.
\end{prop}
\begin{proof}
The proof is similar to the proof of Proposition $3.3$ in \cite{DFL21}, so we will only sketch the argument, for the sake of completeness. 

We note that for each primitive character $\chi$ modulo $q$, we have one of the following three cases: 
\begin{enumerate}
\item $\chi \notin \mathcal{T}_0$,
\item $\chi \in \mathcal{T}_r$, for each $r \leq K$,
\item There exists some $j <K$ such that $\chi \in \mathcal{T}_r$ for each $r \leq j$, but $\chi \notin \mathcal{T}_{j+1}$.
\end{enumerate}
 Note that the first case above is the first condition in the statement of the Proposition. If the second case is satisfied, then we use Lemma \ref{lem_log} with $x=T^{\beta_K}$ and inequality \eqref{ineq_exp}. Finally, if the third condition is satisfied, there exists $u >j$ 
 such that $|\Re P_{I_{j+1}}(\chi;u)| > \ell_{j+1}/(ve^2)$. We then use Lemma \ref{lem_log} with $x=T^{\beta_j}$, and the conclusion follows.
\end{proof}

We also need the following two key lemmas.
\begin{lem}
\label{technical_lemma}
Suppose that $(v+2) \sum_{r=0}^K \ell_r \beta_r <1$ and that for any $j \leq K-1$,  we have $(v+2) \sum_{r=0}^j \ell_r \beta_r + v \sum_{r=j+1}^{K} \ell_r \beta_r + 2s_{j+1} \beta_{j+1} <1$. Let $0 \leq u \leq K$ in $(1)$ and $j+1 \leq u \leq K$ in $(3)$. Suppose that $\beta_{r+1} = e \beta_r$, for $r=0,\ldots, K-1$. Then for any $\varepsilon >0$, we have 

\begin{align*}
(1) & \sumstar_{\chi \pmod q} |  M(\chi ) |^{2v} \Big( \Re P_{I_0}(\chi ;u)\Big)^{2s_0} \ll q e^{v^2K(1+\varepsilon)} \frac{ (2s_0)!v^{2s_0}\sqrt{s_0}}{(s_0)!2^{s_0}} (\log \log q^{\beta_0})^{2s_0} ,\\ 
(2)& \sumstar_{\chi \pmod q} D_{K,v}(\chi ) ^2 |M(\chi )|^{2v} \ll q \Big(  1 + \frac{1}{2^{\ell_0}} (\log q^{\beta_0})^{16v^2}\Big), \\
 (3) &\sumstar_{\chi \pmod q} D_{j,v}(\chi ) ^2 |M(\chi )|^{2v} (\Re P_{I_{j+1}}(\chi ;u))^{2s_{j+1}}  \ll q e^{v^2(K-j)(1+\varepsilon)} \frac{(2s_{j+1})! v^{2s_{j+1}} \sqrt{s_{j+1}}}{(s_{j+1})! 2^{s_{j+1}}} \Big(\sum_{p \in I_{j+1}} \frac{1}{p} \Big)^{2s_{j+1}}\\
 & \qquad \qquad \qquad \quad \quad \times \Big(  1 + \frac{1}{2^{\ell_0}} (\log q^{\beta_0})^{16v^2}\Big).
\end{align*}
\end{lem}
We also prove the following bound on the contribution of primes square.
\begin{lem}
\label{lemma_squares}
Let $S_{j,v}$ be defined as in \eqref{sjk}. For $0 \leq j \leq K$, we have
$$\sumstar_{\chi \pmod q} S_{j,v}(\chi)^2 \ll q.$$
\end{lem} 

We will now first prove Theorem \ref{mol_ub}, and then proceed to the proofs of the technical lemmas above.

\subsection{Choice of parameters and proof of Theorem \ref{mol_ub}}

We choose 
 \begin{align}
 \label{beta_j}
 \beta_j = \frac{e^j \log \log q}{\log q},
 \end{align}
 for $j \leq K$.  We choose $\beta_K=c$, where $c$ is a small positive constant such that
 \begin{equation}
 \label{condition_c}
  c< \frac{4 (e^{1/4}-1)^4}{e(v+2)^4}.
  \end{equation}
 We further choose for $j \leq K$, 
 \begin{equation}
 s_j = 2 \Big[ \frac{1}{8 \beta_j},\Big], \, \ell_j=[s_j^{3/4}].
 \label{sj}
 \end{equation} 
Note that condition \eqref{condition_c} ensures that the conditions in Lemma \ref{technical_lemma} are satisfied.

Using Proposition \ref{main_prop},  we write
\begin{align}
& \sumstar_{\chi \pmod q}  \Big| L (1/2,\chi ) M(\chi ) \Big|^v \leq \sumstar_{\substack{\chi \pmod q \\ \chi \notin \mathcal{T}_0}}  \Big| L (1/2,\chi ) M(\chi ) \Big|^v  \label{bound1} \\
& + \sumstar_{\chi \pmod q }  \exp \Big(  \frac{v (1+\lambda) }{2 \beta_K } \Big) D_{K,v}(\chi ) S_{K,v}(\chi ) |M(\chi )|^v \label{bound2} \\
&+ \sumstar_{\chi \pmod q}  \sum_{\substack{0 \leq j \leq K-1 \\ j < u \leq K}} \exp \Big( \frac{ v(1+\lambda) }{2\beta_j } \Big) D_{j,v} (\chi) S_{j,v}(\chi) \Big(\frac{ve^2 \Re P_{I_{j+1}}(\chi ;u)}{\ell_{j+1}}  \Big)^{s_{j+1}} |M(\chi )|^v. \label{bound3}
\end{align}
For the first term \eqref{bound1}, we have that for some $0 \leq u \leq K$ ,
\begin{align}
&\sumstar_{\substack{\chi \pmod q \\ \chi \notin \mathcal{T}_0}}  \Big| L (1/2,\chi ) M(\chi ) \Big|^v \leq \sumstar_{\chi \pmod q}  \Big| L (1/2,\chi ) M(\chi ) \Big|^v \Big( \frac{ve^2 \Re P_{I_0}(\chi ;u)}{\ell_0} \Big)^{s_0} \nonumber  \\
& \leq \Big( \sumstar_{\chi \pmod q} \Big|L(1/2,\chi ) \Big|^{2v}\Big)^{1/2}  \Big( \sumstar_{\chi \pmod q} |  M(\chi ) |^{2v} \Big( \frac{ve^2 \Re P_{I_0}(\chi ;u)}{\ell_0} \Big)^{2s_0} \Big)^{1/2}. \label{not_to}
\end{align}
Under GRH,  we have the sharp upper bound of Soundararajan \cite{Sou09} and Harper \cite{Har13} for Dirichlet $L$-functions. 
Applying \cite[Theorem 1]{Sza24} (in the notation therein set $2k=2$, $a_1=a_2=v$, $t_1=t_2=0$) we have that
$$ \sumstar_{\chi \pmod q} \Big|L(1/2,\chi ) \Big|^{2v} \ll q (\log q)^{v^2}.$$
Using equation \eqref{not_to}, Lemma \ref{technical_lemma}, the bound above and Stirling's formula, we have
\begin{align*}
\sumstar_{\substack{\chi \pmod q \\ \chi \notin \mathcal{T}_0}} & \Big| L (1/2,\chi ) M(\chi ) \Big|^v \leq q \exp \Big( - \frac{1}{4} s_0  \log s_0 + s_0 \log (2^{1/2}v^2 e^{3/2})+ \frac{1}{4} \log s_0 \\
&+s_0 \log \log \log q^{\beta_0}+\frac{v^2 K(1+\varepsilon)}{2} \Big) (\log q)^{v^2/2} = o(q),
\end{align*} 
where we used the choices \eqref{beta_j} and \eqref{sj}, as well as the fact that $K \asymp \log \log q$.

For the term \eqref{bound2}, we also have that 
\begin{align*}
 \sumstar_{\chi \pmod q}  & D_{K,v}(\chi ) S_{K,v}(\chi ) |M(\chi ) |^v \leq \Big( \sumstar_{\chi \pmod q} D_{K,v}(\chi ) ^2 |M(\chi )|^{2v} \Big)^{1/2} \\ 
& \times  \Big( \sumstar_{\chi \pmod q} S_{K,v}(\chi )^2 \Big)^{1/2} ,
\end{align*}
and for the term \eqref{bound3}, we have
\begin{align}
& \sumstar_{\chi \pmod q}  D_{j,v} (\chi) S_{j,v}(\chi ) \Big(\frac{ve^2 \Re P_{I_{j+1}}(\chi ;u)}{\ell_{j+1}}  \Big)^{s_{j+1}} |M(\chi ) |^v  \label{sumj} \\
& \leq \Big(  \frac{ve^2}{\ell_{j+1}} \Big)^{s_{j+1}}  \Big( \sumstar_{\chi \pmod q} D_{j,v}(\chi ) ^2 |M(\chi )|^{2v}  (\Re P_{I_{j+1}}(\chi ;u))^{2s_{j+1}}\Big)^{1/2}  \Big( \sumstar_{\chi \pmod q} S_{j,v}(\chi )^2 \Big)^{1/2}. \nonumber 
\end{align}

We focus on bounding the term \eqref{bound3}; bounding \eqref{bound2} follows similarly. Using \eqref{sumj}, Lemma \ref{technical_lemma} and Stirling's formula, and since 
$$  \sum_{p \in I_{j+1}} \frac{1}{p} = 1+ O \Big(  \frac{1}{\beta_j \log q}\Big) \leq 2,$$
we have that
\begin{align}
& \eqref{bound3}   \ll q  \sum_{\substack{0 \leq j \leq K-1 }} (K-j) \exp \Big( \frac{ v(1+\lambda) }{2\beta_j } \Big) \exp \Big( -\frac{1}{4}s_{j+1} \log s_{j+1} +s_{j+1} \log (2^{3/2} v^2 e^{3/2})  \nonumber  \\
&+\frac{1}{4}\log s_{j+1}+ \frac{v^2(K-j)(1+\varepsilon)}{2} \Big) \ll q \sum_{0 \leq j <K} \exp \Big( \frac{A \log q}{e^j \log \log q} - \frac{\log q}{16e^{j+1}}+ \frac{(j+1) \log q}{16e^{j+1} \log \log q} \nonumber \\
&+ \frac{ (\log q) (\log \log \log q)}{16e^{j+1} \log \log q}  \frac{ \log \log q -(j+1)- \log \log \log q}{4}+ \frac{v^2(K-j)(1+\varepsilon)}{2}+\log(K-j) \Big),
\nonumber 
\end{align}
where $A = \frac{v(1+\lambda)}{2} + \frac{\log 4}{16 e} + \frac{ \log (2^{3/2} v^2 e^{3/2})}{4e}.$
Rearranging the above, we have
\begin{align*}
\eqref{bound3} & \ll q \sum_{0 \leq j <K} \exp \Big(  \frac{1}{16e^{j+1}} \Big( \frac{(j+1) \log q}{\log \log q} - \log q + \frac{16Ae \log q}{\log\log q}+ \frac{( \log q )(\log \log \log q)}{\log \log q} \Big) \\
& + \frac{ \log \log q -(j+1)- \log \log \log q}{4} + \frac{v^2(K-j)(1+\varepsilon)}{2} + \log (K-j)\Big).
\end{align*}

Now using the fact that $K< \log \log q - \log \log \log q + \log d $, and for $d$ small enough, note that we have
\begin{align*}
\eqref{bound3} &\ll q \sum_{0 \leq j <K} \exp \Big(  - \frac{b \log q}{e^{j+1} \log \log q}+ \frac{ \log \log q -(j+1)- \log \log \log q}{4} + \frac{v^2(K-j)(1+\varepsilon)}{2} \\
&+ \log (K-j)\Big),
\end{align*}
for some $b>0$. Bounding $\log(K-j) \leq (K-j)$, we get that the above is
\begin{align*}
\eqref{bound3}& \ll q \sum_{0 \leq j <K}  \exp \Big(  - \frac{b \log q}{e^{j+1} \log \log q}+ \frac{ \log \log q -(j+1)- \log \log \log q}{4} + (K-j) \Big( \frac{v^2(1+\varepsilon)}{2}+1\Big) \Big) \\
& \ll q \exp \Big( \frac{ \log \log q - \log \log \log q}{4} + K \Big( \frac{v^2(1+\varepsilon)}{2}+1\Big) \Big) \\
& \times \int_{0}^{K-1} \exp \Big(- \frac{b \log q}{e^{t+1} \log \log q} - \frac{t}{4} - t \Big(\frac{v^2(1+\varepsilon)}{2}+1\Big) \Big) \, dt.
\end{align*}
With the change of variables $e^{-t} = \frac{ e \log \log q }{b \log q} u$, it follows that
\begin{align*}
\eqref{bound3} & \ll  q \exp \Big( \frac{ \log \log q - \log \log \log q}{4} + K \Big( \frac{v^2(1+\varepsilon)}{2}+1\Big) \Big) \Big( \frac{ \log \log q}{\log q} \Big)^{ \frac{1}{4} + \frac{v^2(1+\varepsilon)}{2}+1} \\
& \times  \int_{ \frac{ b \log q}{e^K \log \log q}}^{\frac{b \log q}{e \log \log q}} \exp(-u) u^{ \frac{1}{4} + \frac{k^2(1+\varepsilon)}{2}} \, du \\
& \leq q \exp \Big( \frac{ \log \log q - \log \log \log q}{4} + K \Big( \frac{v^2(1+\varepsilon)}{2}+1\Big) \Big) \Big( \frac{ \log \log q}{\log q} \Big)^{ \frac{1}{4} + \frac{v^2(1+\varepsilon)}{2}+1}  \Gamma \Big( \frac{5}{4}+ \frac{v^2(1+\varepsilon)}{2} \Big).
\end{align*}
Since $e^K \asymp \log q/\log \log q$, from the above, we get that
\begin{align*}
\eqref{bound3} \ll q.
\end{align*}

\kommentar{$$ \frac{\log \log q}{4} - \frac{j+1}{4} <  \frac{ \varepsilon \log q}{16e^{j+1}},$$ and 
$$   \frac{A \log q}{e^j \log \log q}  + \frac{ (\log q) (\log \log \log q)}{16e^{j+1} \log \log q} - \frac{\log \log \log q}{4} + \frac{k^2(K-j)(1+\varepsilon)}{2}+\log(K-j) <  \frac{ \varepsilon \log q}{16e^{j+1}} .$$
After a relabeling of $\varepsilon$, it follows that 
\begin{align}
\eqref{bdint} \ll q \sum_{0 \leq j <K} \exp \Big(\frac{-(1-\varepsilon) \log q}{16 e^{j+1}} + \frac{(j+1) \log q}{16e^{j+1}\log \log q}  \Big) 
\end{align}
}

Similarly 
$$ \eqref{bound2} \ll q.$$
This finishes the proof of Theorem \ref{mol_ub}.
\qed

\subsection{Proof of Lemma \ref{technical_lemma}} We are now ready to begin the proof of the technical Lemma \ref{technical_lemma}.
We start with the second part of the lemma. The term we need to bound is
\begin{align}
\label{to_bd}
\sumstar_{\chi \pmod q} &  \prod_{r=0}^K (1+e^{-\ell_r/2})^2 \mathcal{E}_r (\chi ),
\end{align}
where
\begin{align}
\mathcal{E}_r(\chi ) &= \sum_{\substack{ p | m_{r_1} n_{r_1} m_{r_2} n_{r_2} \Rightarrow p \in I_r \\ \Omega( m_{r_1} n_{r_1}) \leq \ell_r \\ \Omega(m_{r_2} n_{r_2}) \leq \ell_r}} \frac{ (v/2)^{\Omega(m_{r_1}n_{r_1}m_{r_2} n_{r_2})} \chi (n_{r_1} n_{r_2} ) \overline{ \chi}(m_{r_1} m_{r_2}) }{\sqrt{m_{r_1} m_{r_2} n_{r_1} n_{r_2}}} \nonumber  \\
& \times a(m_{r_1} m_{r_2}n_{r_1} n_{r_2};K)  \nu(n_{r_1}) \nu(n_{r_2}) \nu(m_{r_1}) \nu(m_{r_2}) \nonumber \\
& \times \sum_{\substack{p|f_rh_r \Rightarrow p \in I_r \\ \Omega(f_r) \leq v \ell_r \\ \Omega(h_r) \leq v \ell_r}} \frac{\chi(f_r) \overline{\chi}(h_r) \lambda(f_r) \lambda(h_r) a(f_r h_r;K) \nu_v(f_r;\ell_r) \nu_v(h_r;\ell_r) }{\sqrt{f_rh_r}}. \label{epsilonr}
\end{align}
Note that the summands in \eqref{to_bd} are positive, so we bound \eqref{to_bd} by the sum over all characters $\chi$ modulo $q$. Using orthogonality of characters, we only keep the term with $ n_{r_1} n_{r_2} f_r \equiv   m_{r_1} m_{r_2} h_r \pmod q$. Since $\prod_{r=0}^K n_{r_1} n_{r_2} f_r <q$ and $\prod_{r=0}^K m_{r_1} m_{r_2} h_r <q$ (because $ (v+2) \sum_{r=0}^K \beta_r \ell_r <1$), this can happen only when $\prod_{r=0}^K n_{r_1} n_{r_2} f_r=\prod_{r=0}^K m_{r_1} m_{r_2} h_r$, hence when $n_{r_1} n_{r_2} f_r= m_{r_1} m_{r_2} h_r$ for each $r \leq K$. 

Then the term we need to bound is 
\begin{align*}
q \prod_{r=0}^K &E_r := q \prod_{r=0}^K  \sum_{\substack{ p | m_{r_1} n_{r_1} m_{r_2} n_{r_2}f_rh_r  \Rightarrow p \in I_r \\ \Omega( m_{r_1} n_{r_1}) \leq \ell_r \\ \Omega(m_{r_2} n_{r_2}) \leq \ell_r \\\Omega(f_r) \leq v \ell_r \\ \Omega(h_r) \leq v \ell_r \\ n_{r_1} n_{r_2} f_r = m_{r_1} m_{r_2} h_r }} \frac{ (v/2)^{\Omega(m_{r_1}n_{r_1}m_{r_2} n_{r_2})}a(m_{r_1} m_{r_2}n_{r_1} n_{r_2};K)   }{\sqrt{m_{r_1} m_{r_2} n_{r_1} n_{r_2} f_rh_r}} \\
& \times  \nu(n_{r_1}) \nu(n_{r_2}) \nu(m_{r_1}) \nu(m_{r_2}) \lambda(f_r) \lambda(h_r)a(f_r h_r;K)  \nu_v(f_r;\ell_r) \nu_v(h_r;\ell_r) .
\end{align*}

If $\max\{ \Omega(m_{r_1} n_{r_1}), \Omega(m_{r_2} n_{r_2})\} >\ell_r$ or if $\max\{\Omega(f_r),\Omega(h_r)\} >v \ell_r$ then, since $v \geq 1$, we have that $2^{\Omega(m_{r_1} m_{r_2} n_{r_1} n_{r_2} f_r h_r)} \geq 2^{\ell_r}$.  It then follows that 
\begin{align*}
E_r & =  \sum_{\substack{ p | m_{r_1} n_{r_1} m_{r_2} n_{r_2}f_rh_r  \Rightarrow p \in I_r \\ n_{r_1} n_{r_2} f_r= m_{r_1} m_{r_2} h_r  }} \frac{ (v/2)^{\Omega(m_{r_1}n_{r_1}m_{r_2} n_{r_2})}a(m_{r_1} m_{r_2} n_{r_1} n_{r_2};K)  }{\sqrt{m_{r_1} m_{r_2} n_{r_1} n_{r_2} f_rh_r}} \\
& \times  \nu(n_{r_1}) \nu(n_{r_2}) \nu(m_{r_1}) \nu(m_{r_2})\lambda(f_r) \lambda(h_r)  a(f_r h_r;K)  \nu_v(f_r) \nu_v(h_r)  \\
&+ O \Bigg(  \frac{1}{2^{\ell_r}}  \sum_{\substack{ p | m_{r_1} n_{r_1} m_{r_2} n_{r_2}f_rh_r  \Rightarrow p \in I_r \\ n_{r_1} n_{r_2} f_r = m_{r_1} m_{r_2} h_r }} \frac{2^{\Omega(m_{r_1} m_{r_2} n_{r_1} n_{r_2} f_rh_r)} (v/2)^{\Omega(m_{r_1}n_{r_1}m_{r_2} n_{r_2})} }{\sqrt{m_{r_1} m_{r_2} n_{r_1} n_{r_2} f_rh_r}}   \nu(n_{r_1}) \nu(n_{r_2}) \nu(m_{r_1}) \nu(m_{r_2})\\
& \times \nu_v(f_r) \nu_v(h_r) \Bigg),
\end{align*}
where we used the fact that $a(n;K) \leq 1$ for any $n$.
 Now let $M_r=m_{r_1} m_{r_2}$ and $N_r= n_{r_1} n_{r_2}$.  We then rewrite
 \begin{align}
 E_r &= \sum_{\substack{p|M_rN_r f_rh_r \Rightarrow p \in I_r \\ N_rf_r=M_rh_r}}  \frac{ (v/2)^{\Omega(M_rN_r)}a(M_r N_r;K)   }{\sqrt{M_rN_r f_rh_r}} \nu_2(N_r) \nu_2(M_r)  \lambda(f_r) \lambda(h_r) a(f_r h_r;K)  \nonumber  \\
& \times   \nu_v(f_r) \nu_v(h_r) + O \Bigg(  \frac{1}{2^{\ell_r}}  \sum_{\substack{ p | M_rN_rf_rh_r  \Rightarrow p \in I_r \\ N_r f_r = M_r h_r }} \frac{ (2v)^{\Omega( M_rN_rf_rh_r)} }{\sqrt{M_rN_r f_rh_r}}   \Bigg), \label{er}
 \end{align}
where we used the fact that $\nu_2(f) \leq 2^{\Omega(f)}$ and $\nu_v(f) \leq v^{\Omega(f)}$.

Let $A_r=(N_r, M_r)$ and $B_r=(f_r,h_r)$, and we write $N_r= A_r C_r$,  and $M_r=A_rD_r$ with $(C_r,D_r)=1$.  Similarly write $f_r=B_rf_{r1}$ and $h_r=B_rh_{r1}$ with $(f_{r1}, h_{r1})=1$.  Since $N_rf_r=M_rh_r$,  it follows that $f_{r1}=D_r$ and $h_{r1}=C_r$.

The error term in the evaluation of $E_r$ in equation \eqref{er} then becomes 
\begin{align}
\frac{1}{2^{\ell_r}} &\sum_{\substack{p| A_rB_rC_rD_r \Rightarrow p \in I_r \\ (C_r,D_r)=1}} \frac{ (2v)^{\Omega(A_r^2B_r^2C_r^2 D_r^2 )}}{A_rB_rC_rD_r} \leq \frac{1}{2^{\ell_r}} \Big( \sum_{p|A \Rightarrow p \in I_r}\frac{(2v)^{\Omega(A^2)}}{A} \Big)^4 \nonumber  \\
& \leq \frac{1}{2^{\ell_r}}\prod_{p \in I_r} \Big(1 - \frac{4v^2}{p} \Big)^{-4} =
\begin{cases}
O \Big( \frac{1}{2^{\ell_r}} \Big) & \mbox{ if } r \neq 0 \\
O \Big(  \frac{1}{2^{\ell_0}} (\log q^{\beta_0})^{16v^2}\Big) & \mbox{ if } r=0.
\end{cases}
\label{error_er}
\end{align}

We now focus on the main term in \eqref{er}.  We rewrite it as
\begin{align*}
 \sum_{\substack{p| A_rB_rC_rD_r \Rightarrow p \in I_r \\ (C_r,D_r)=1}}  &  \frac{ (v/2)^{\Omega(A^2 C_rD_r)} a(A_r^2B_r^2C_r^2D_r^2 ;K)}{A_rB_rC_rD_r} \nu_2(A_rC_r) \nu_2(A_rD_r)\lambda(C_rD_r) \\
 & \times \nu_v(B_rD_r) \nu_v(B_rC_r) = \prod_{p \in I_r} \Big( 1 + O \Big(\frac{1}{p^2}\Big) \Big).
\end{align*}

Combining the above and equation \eqref{error_er}, it follows that
\begin{equation*}
q \prod_{r=0}^K E_r  \ll q \Big( 1 +\frac{1}{2^{\ell_0}} (\log q^{\beta_0})^{16v^2} \Big),
\end{equation*} which finishes the proof of the second bound in the lemma.

We now focus on the third bound.  We need to bound
\begin{align*}
\sumstar_{\chi \pmod q } \prod_{r=0}^j (1+e^{-\ell_r/2})^2 \mathcal{E}_r(\chi) \prod_{r=j+1}^K \mathcal{E}_r(\chi),
\end{align*}
where if $r \leq j$,  then $\mathcal{E}_r(\chi)$ is given as in the previous case by equation \eqref{epsilonr}.

 If $r=j+1$,  then
\begin{align*}
\mathcal{E}_{r}(\chi) &= \frac{ (2s_{j+1})!}{4^{s_{j+1}}}  \sum_{\substack{p|f_rh_r c_rd_r\Rightarrow p \in I_r \\ \Omega(f_r) \leq v \ell_r \\ \Omega(h_r) \leq v \ell_r \\ \Omega(c_rd_r)=2s_{j+1}}}  \frac{\chi(f_r) \overline{\chi}(h_r) \lambda(f_r) \lambda(h_r) a(f_r h_r;K) \nu_v(f_r;\ell_r) \nu_v(h_r;\ell_r) }{\sqrt{f_rh_r}} \\
& \times \chi(c_r) \overline{\chi(d_r)} a(c_rd_r;u)   \nu(c_r) \nu(d_r).
\end{align*}

If $r \geq j+2$,  then 
\begin{align*}
\mathcal{E}_{r}(\chi) &=  \sum_{\substack{p|f_rh_r \Rightarrow p \in I_r \\ \Omega(f_r) \leq v \ell_r \\ \Omega(h_r) \leq v \ell_r }}  \frac{\chi(f_r) \overline{\chi}(h_r) \lambda(f_r) \lambda(h_r) a(f_r h_r;K) \nu_v(f_r;\ell_r) \nu_v(h_r;\ell_r) }{\sqrt{f_rh_r}}.
\end{align*}
Since all the summands are positive, we extend the sum over primitive characters to a sum over all $\chi \pmod q$,  and use orthogonality of characters. Since
\begin{equation}
(v+2) \sum_{r \leq j} \ell_r \beta_r + 2s_{j+1} \beta_{j+1}+ \sum_{r=j+1}^K v \ell_r \beta_r <1.
\label{cond_sum}
\end{equation}
it follows that we have to bound the following expression
\begin{align*}
q \prod_{r =0}^K E_r,
\end{align*}
where if $r \leq j$, as before, we have that
\begin{align*}
E_r &= \sum_{\substack{ p | m_{r_1} n_{r_1} m_{r_2} n_{r_2}f_rh_r  \Rightarrow p \in I_r \\ \Omega( m_{r_1} n_{r_1}) \leq \ell_r \\ \Omega(m_{r_2} n_{r_2}) \leq \ell_r \\\Omega(f_r) \leq v \ell_r \\ \Omega(h_r) \leq v \ell_r \\ n_{r_1} n_{r_2} f_r = m_{r_1} m_{r_2} h_r }} \frac{ (k/2)^{\Omega(m_{r_1}n_{r_1}m_{r_2} n_{r_2})}a(m_{r_1} m_{r_2}n_{r_1} n_{r_2};K)   }{\sqrt{m_{r_1} m_{r_2} n_{r_1} n_{r_2} f_rh_r}} \\
& \times  \nu(n_{r_1}) \nu(n_{r_2}) \nu(m_{r_1}) \nu(m_{r_2}) \lambda(f_r) \lambda(h_r)a(f_r h_r;K)  \nu_v(f_r;\ell_r) \nu_v(h_r;\ell_r) .
\end{align*}
When $r=j+1$,  
\begin{align*}
E_r &= \frac{ (2s_{r})!}{4^{s_{r}}}  \sum_{\substack{p|f_rh_r c_rd_r\Rightarrow p \in I_r \\ \Omega(f_r) \leq v \ell_r \\ \Omega(h_r) \leq v \ell_r \\ \Omega(c_rd_r)=2s_{r}\\ f_rc_r=h_rd_r}}  \frac{  \lambda(f_r) \lambda(h_r) a(f_r h_r;K) \nu_v(f_r;\ell_r) \nu_v(h_r;\ell_r) a(c_rd_r;u)   \nu(c_r) \nu(d_r)}{\sqrt{f_rh_rc_rd_r}}  ,
\end{align*}
and when $r \geq j+2$,
\begin{equation}
\label{bigr}
E_r =  \sum_{\substack{p|f_rh_r \Rightarrow p \in I_r \\ \Omega(f_r) \leq v \ell_r \\ \Omega(h_r) \leq v \ell_r \\ f_r=h_r}} \frac{\lambda(f_r) \lambda(h_r) a(f_r h_r;K) \nu_v(f_r;\ell_r) \nu_v(h_r;\ell_r) }{\sqrt{f_rh_r}}.
\end{equation}

We deal with the product $\prod_{r=0}^j E_r$ in exactly the same way as before, and it follows that
\begin{equation}
\prod_{r=0}^j E_r \ll 1 +\frac{1}{2^{\ell_0}} (\log q^{\beta_0})^{16v^2}.
\label{er1}
\end{equation}
To bound $E_{r}$ when $r=j+1$,  we get
$$E_r \leq  \frac{ (2s_{r})!}{4^{s_{r}}}  \sum_{\substack{p|f_rh_r c_rd_r\Rightarrow p \in I_r \\ \Omega(f_r) \leq v \ell_r \\ \Omega(h_r) \leq v \ell_r \\ \Omega(c_rd_r)=2s_{r}\\ f_rc_r=h_rd_r}} \frac{ \nu_v(f_r) \nu_v(h_r) \nu(c_r) \nu(d_r)}{\sqrt{f_rh_rc_rd_r}}  .$$
Using a similar change of variables as before,  we write $f_r=B_rf_{r1}, c_r = A_rh_{r1}, h_r=B_rh_{r1}, d_r=A_rf_{r1}$ with $(f_{r1}, h_{r1})=1$.  Then we have
\begin{align*}
E_r \leq \frac{ (2s_{r})!}{4^{s_{r}}}  \sum_{\substack{p |A_rB_rf_{r1}h_{r1} \Rightarrow p \in I_r \\ \Omega(A_rf_{r1}) \leq v\ell_r \\ \Omega(B_rh_{r1})\leq v \ell_r \\ \Omega(A_r^2f_{r1} h_{r1})=2s_r}} \frac{\nu_v(B_rf_{r1}) \nu_v(B_rh_{r1})\nu(A_rf_{r1}) \nu(A_rh_{r1})}{A_rB_rf_{r1}h_{r1}}.
\end{align*}
Now we use the fact that 
$$\nu_v(B_rf_{r1}) \leq d_v(B_rf_{r1}) \leq d_v(B_r) d_v(f_{r1}) \leq d_v(B_r) v^{\Omega(f_{r1})},$$ and a similar inequality holds for $\nu_v(B_rh_{r1})$.  We also use the fact that $\nu(A_rf_{r1}) \leq \nu(A_r) \nu(f_{r1})$ and similarly for $\nu(A_rh_{r1})$. We then get that
\begin{align}
E_r \leq \frac{ (2s_{r})!}{4^{s_{r}}}  \sum_{\substack{p |A_rB_rf_{r1}h_{r1} \Rightarrow p \in I_r \\ \Omega(A_r^2f_{r1} h_{r1})=2s_r}} \frac{d_v(B_r)^2 v^{2s_r-2\Omega(A_r)} \nu(A_r) \nu(f_{r1}) \nu(h_{r1})}{A_rB_rf_{r1}h_{r1}}, \label{er4}
\end{align}
where we also used the fact that $\nu(A_r)^2 \leq \nu(A_r)$.
We first consider the sum over $f_{r1}$.  We have 
\begin{align}
\sum_{\substack{p|f_{r1} \Rightarrow p \in I_r \\ \Omega(f_{r1}) = 2s_r-\Omega(h_{r1})-2\Omega(A_r)}}  \frac{ \nu(f_{r1})}{f_{r1}}  = \frac{1 }{(2s_r-\Omega(h_{r1})-2\Omega(A_r))!}\Big(  \sum_{p \in I_r} \frac{1}{p} \Big)^{2s_r-\Omega(h_{r1}) - 2\Omega(A_r)}.
\label{sumfr1}
\end{align}
Let $S= \sum_{p \in  I_r} \frac{1}{p}.$ We introduce the sum over $h_{r1}$,  and we have 
\begin{align}
\sum_{\substack{p|h_{r1} \Rightarrow p \in I_r\\\Omega(h_{r1}) \leq 2s_r-2\Omega(A_r)}} &\frac{ \nu(h_{r1})}{(2s_r-\Omega(h_{r1})-2\Omega(A)!h_{r1}S^{\Omega(h_{r1})}} = \sum_{i=0}^{2s_r-2\Omega(A_r)} \frac{1}{(2s_r-2\Omega(A)-i)!i!S^i}  \Big( \sum_{p \in I_r} \frac{1}{p} \Big)^i \nonumber  \\
& = \frac{2^{2s_r-2\Omega(A_r)}}{(2s_r-2\Omega(A_r))!}. \label{sumhr1}
\end{align}
Now we introduce the sum over $A_r$ in \eqref{er4}. Using the above equations,  it follows that
\begin{align*}
\sum_{\substack{p|A_r \Rightarrow p \in I_r \\ \Omega(A_r) \leq s_r}} & \frac{\nu(A_r)}{ (2s_r-2\Omega(A_r))!(2vS)^{2\Omega(A)}A_r} = \sum_{i=0}^{s_r} \frac{1}{ (2s_r-2i)!i!(2vS)^{2i}} S^i   .
\end{align*}
Using Stirling's formula, we get that the above is
\begin{align}
\label{sumar}
\ll \sum_{i=0}^{s_r} \frac{ \sqrt{s_r-i}}{4^{2s_r-2i} i!(s_r-i)!^2 (4v^2S)^i} \ll \frac{\sqrt{s_r}}{4^{2s_r}(s_r)!} \Big( 1+ \frac{4}{v^2S} \Big)^{s_r} \ll \frac{\sqrt{s_r}}{(s_r)! 2^{s_r}},
\end{align}
where we used the fact that $v \geq 1$ and that $S=1+o(1)$ (which follows from the fact that $\beta_{r+1}= e \beta_r$).
We finally introduce the sum over $B_r$ and we have
\begin{equation}
\label{sumbr}
\sum_{p|B_r \Rightarrow p \in I_r} \frac{d_v(B_r)^2}{B_r} \ll \prod_{p \in I_r} \Big(1+\frac{1}{p} \Big)^{v^2} =O(1).
\end{equation}
From equations \eqref{er4}, \eqref{sumfr1}, \eqref{sumhr1}, \eqref{sumar}, \eqref{sumbr} it follows that when $r=j+1$,
\begin{equation}
E_r \leq \frac{(2s_r)! (vS)^{2s_r}\sqrt{s_r}}{(s_r!) 2^{s_r}}.
\label{erj+1}
\end{equation}

Now we consider the contribution from $r \geq j+2$ in equation \eqref{bigr}.  In this case,  we use the fact that $\nu_v(f_r;\ell_r) \leq \nu_v(f_r) \leq v^{\Omega(f_r)}$, and we get that
\begin{align*}
E_r \leq \sum_{p|f_r \Rightarrow p \in I_r} \frac{v^{2\Omega(f_r)}}{f_r} \ll \prod_{p \in I_r} \Big(1+ \frac{v^2}{p}\Big) \leq e^{v^2(1+\varepsilon)},
\end{align*} 
hence
\begin{equation} 
\label{erj+2}
\prod_{r=j+2}^K E_r \ll e^{v^2(K-j)(1+\varepsilon)}.
\end{equation}
Using equations \eqref{er1}, \eqref{erj+1} and \eqref{erj+2}, the conclusion follows for the third bound of the lemma. 

Bounding $(1)$ in Lemma \ref{technical_lemma} is similar, except that in equation \eqref{erj+1} we bound 
$$\Big( \sum_{p \in I_0} \frac{1}{p} \Big)^{2s_0} \ll (\log \log q^{\beta_0})^{2s_0}.$$
\qed


\subsection{Proof of Lemma \ref{lemma_squares}}
The proof is similar to the proof of Lemma $5.1$ in \cite{DFL21}.  Let
$$F_m(\chi;j) = \sum_{2^m<p\leq 2^{m+1}} \frac{\chi(p^2)b(p;j)}{p},$$
where recall that $b(p;j)$ is given in equation \eqref{bp}. 
For simplicity,  we will denote $F_m(\chi;j)$ by $F_m(\chi)$.  Let
$$\mathcal{F}(m)  = \Big\{  \chi :  |\Re F_m(\chi) |> \frac{1}{2^{m/4}} \text{ but } |\Re F_n(\chi)| \leq \frac{1}{2^{n/4}}, \forall m+1 \leq n \leq \frac{\beta_j \log q}{2 \log 2} \Big\}.$$
Note that the sets $\mathcal{F}(m)$ are disjoint,  and we have
\begin{equation}
 \sumstar_{\chi \pmod q} S_{j,v}(\chi)^2 \leq \sum_{m \leq \frac{\beta_j \log q}{2\log 2}} \sum_{\chi \in \mathcal{F}(m)} S_{j,v}(\chi)^2 + \sum_{\chi \notin \mathcal{F}(m),  \forall m} S_{j,v}(\chi)^2.
 \label{sumsquares}
 \end{equation}
If $\chi \notin \mathcal{F}(m)$ for any $m$, then note that $|\Re F_n(\chi)| \leq \frac{1}{2^{n/4}}$ for all $n \leq \frac{\beta_j \log q}{2 \log 2}$, and hence in this case, 
$$S_{j,v}(\chi)^2 = O(1).$$
If $\chi \in \mathcal{F}(m)$ for some $m \leq \frac{\beta_j \log q}{2 \log 2}$, then 
$$\sum_{\ell=0}^{\frac{\beta_j \log q}{2 \log 2}} \Re F_{\ell}(\chi) \leq \frac{1}{2}+ \sum_{\ell=1}^{m} \frac{1}{\ell}+ \sum_{\ell=m+1}^{\frac{\beta_j \log q}{2\log 2}} \frac{1}{2^{\ell/4}} \leq \frac{1}{2}+\log m+1+ \frac{1}{2^{m/4}(1-2^{-1/4})}.$$
It then follows that
$$S_{j,v}(\chi)^2 \leq \exp \Big(v+2v(\log m+1)+ \frac{2v}{2^{m/4}(1-2^{-1/4})}\Big).$$
If $\chi \in \mathcal{F}(m)$, we also have that $(2^{m/4} \Re F_m(\chi))^4 > 1$, so 
\begin{align}
\label{sum_m}
\sum_{\chi \in \mathcal{F}(m)} S_{j,v}(\chi)^2 \leq \exp \Big(v+2v(\log m+1)+ \frac{2v}{2^{m/4}(1-2^{-1/4})}\Big) \sum_{\chi \pmod q} (2^{m/4} \Re F_m(\chi))^4.
\end{align}
Now we have
\begin{align}
\sum_{\chi \pmod q} (2^{m/4} \Re F_m(\chi))^4 = \frac{4!2^m}{2^4} \sum_{\substack{p|fh \Rightarrow 2^m<p\leq 2^{m+1} \\ \Omega(fh)=4}} \frac{ \chi(f^2) \overline{\chi}(h^2) b(f;j) b(h;j)  \nu(f) \nu(h)}{fh}.
\label{sumchi}
\end{align}
We use orthogonality of characters (because $\beta_j \leq \beta_K$ and $\beta_K$ is small enough) to get that the above is 
\begin{align}
\label{eq4}
\eqref{sumchi} = q\frac{4! 2^m}{2^4} \sum_{\substack{p|f \Rightarrow 2^m<p \leq 2^{m+1} \\ \Omega(f)=2}} \frac{ b(f;j)^2 \nu(f)^2}{f^2} \leq q\frac{4! 2^m}{2^5} \Big( \sum_{2^m<p \leq 2^{m+1}} \frac{1}{p^2}\Big)^2\leq q\frac{4!}{2^{m+5}},
\end{align}
where we used the fact that $|b(f;j)| \leq 1$. 

Now using equations \eqref{sumsquares}, \eqref{sum_m} and \eqref{eq4}, we get that 
\begin{align*}
\sumstar_{\chi \pmod q} S_{j,v}(\chi)^2 & \ll q  \sum_{m \leq \frac{\beta_j \log q}{2 \log 2}} \exp \Big(v+2v(\log m+1)+ \frac{2v}{2^{m/4}(1-2^{-1/4})}\Big) \frac{1}{2^{m}} \\
& \ll q \sum_{m=0}^{\infty} \frac{m^{2v}}{2^m} \ll q,
\end{align*}
where in the last line above we used the fact that for $x<1$, we have
$$ \sum_{m=1}^{\infty} m^v x^m \leq \frac{x v!}{(1-x)^{v+1}}.$$
\qed

\section{Proof of Theorem \ref{mol_lb}}
\label{lb}

Using Theorem \ref{thm:twisted-moment-full-orbit}, we have that

\begin{align}
& \frac{1}{|\co|} \sum_{\substack{\chi\in\co}}L(\tfrac12,\chi\eta_1)L(\tfrac12,\overline{\chi\eta_2})M(\chi\eta_1)M(\overline{\chi\eta_2}) \nonumber \\
&=  \sum_{\substack{ M= \prod_{j=0}^K m_j, \, N= \prod_{j=0}^K n_j \\ p|m_j n_j \Rightarrow p \in I_j \\ \Omega(m_j)\leq \ell_j, \, \Omega(n_j) \leq \ell_j}} \frac{ a(M;K) a(N;K) \lambda(M) \lambda(N) \nu(M) \nu(N)}{\sqrt{MN}}\nonumber\\
&\qquad\qquad\times \frac{1}{|\co|} \sum_{\substack{\chi\in\co}}L(\tfrac12,\chi\eta_1)L(\tfrac12,\overline{\chi\eta_2})\chi\eta_1(M)\overline{\chi\eta_2}(N) \nonumber \\
&= L(1,\eta)  \prod_{j=0}^K \sum_{\substack{ p |g_j m_j n_j \Rightarrow p \in I_j \\ (m_j,n_j)=1 \\ \Omega(g_j m_j) \leq \ell_j, \, \Omega(g_j n_j) \leq \ell_j}}  \frac{a^2(g_j;K) a(m_j; K) a(n_j;K) \lambda(m_j) \lambda(n_j) \nu(g_j m_j) \nu(g_j n_j) \eta(g_j m_j n_j)}{g_j m_j n_j} \nonumber \\
&\qquad\qquad +O\Big( \sum_{\substack{ M= \prod_{j=0}^K m_j, \, N= \prod_{j=0}^K n_j \\ p|m_j n_j \Rightarrow p \in I_j \\ \Omega(m_j)\leq \ell_j, \, \Omega(n_j) \leq \ell_j}} \frac{ a(M;K) a(N;K) }{\sqrt{MN}} E(M,N) \Big), \label{mol_mom}
\end{align}
with $E(M,N)$ as in Theorem \ref{thm:twisted-moment-full-orbit}, and where recall that $\eta = \eta_1 \overline{\eta_2}$. We write 
\begin{align*}
& \frac{1}{|\co|} \sum_{\substack{\chi\in\co}}L(\tfrac12,\chi\eta_1)L(\tfrac12,\overline{\chi\eta_2})M(\chi\eta_1)M(\overline{\chi\eta_2}) \nonumber \\ 
&=: \mathcal{M}+\mathcal{E},
\end{align*}
where $\mathcal{M}$ corresponds to the main term above and where $\mathcal{E}$ comes from the error term $E(M, N)$. Using Theorem \ref{thm:twisted-moment-full-orbit}, and since 
$$ \max\bigg(\frac{\log M}{\log q}, \frac{\log N}{\log q}\bigg) \leq \sum_{j=0}^K \ell_j \beta_j < 2 \beta_K^{1/4},$$
we have that
$$ \mathcal{E} \ll q^{-1/4+2\beta_K^{1/4}+\varepsilon},$$ 
recalling that 
$\beta_K$ is given in equation \eqref{condition_c} (with $v=4$). Since we can choose $c$ small enough, the error term above is negligible.

We now focus on the main term $\mathcal{M}$ in the expression \eqref{mol_mom}. For $j=0,\ldots, K$, let $T(j)$ denote term corresponding to $I_j$ in the product over $j$. 

For $j \geq 1$, we write
$$T(j) = U(j)-V(j),$$ where 
$$U(j) = \sum_{\substack{ p |g_j m_j n_j \Rightarrow p \in I_j \\ (m_j,n_j)=1 }} \frac{a(g_j;K)^2 a(m_j; K) a(n_j;K) \lambda(m_j) \lambda(n_j) \nu(g_j m_j) \nu(g_j n_j)\eta(g_j m_j n_j)}{g_j m_j n_j},$$
 and where $V_j$ corresponds to the sum over $g_j,m_j,n_j$ where at least one of $\Omega(g_j m_j)$ or $\Omega(g_j n_j)$ is greater than $\ell_j$. 
Note that we have 
\begin{align}
\label{uj}
U(j) &= \prod_{p \in I_j} \Bigg( 1 + \frac{a^2(p;K)  \eta(p)}{p} - \frac{2a(p;K) \eta(p)}{p}  + \sum_{\substack{i+j+k \geq 2 \\ jk=0}} \frac{a(p^{2i+j+k};K) (-1)^{j+k}\eta(p^{i+j+k}) }{(i+j)! (i+k)! p^{i+j+k}} \Bigg).
\end{align} 
 We also have that 
 \begin{align*}
 |T(j)| \geq |U(j)| - \frac{1}{2^{\ell_j}} \sum_{p |g_j m_j n_j \Rightarrow p \in I_j} \frac{2^{\Omega(g_j^2 m_j n_j)} }{g_j m_j n_j},
 \end{align*}
where we used the fact that $a(r;K) \leq 1$ and $\nu(r) \leq 1$, for any integer $r$. From the above it follows that
$$|T(j)| \geq |U(j)| - \frac{1}{2^{\ell_j}} \prod_{p \in I_j} \Big(1-\frac{4}{p} \Big)^{-1} \Big( 1 - \frac{2}{p} \Big)^{-2}.$$
Hence we get that
\begin{align}
\label{prod_tj}
\Big| \prod_{j=1}^K T(j)  \Big| \geq \Big| \prod_{j=1}^K U(j) \Big| \prod_{j=1}^K \Bigg(1 - \frac{1}{2^{\ell_j} |U(j)|} \prod_{p \in I_j} \Big(1-\frac{4}{p} \Big)^{-1} \Big( 1 - \frac{2}{p} \Big)^{-2} \Bigg).
\end{align}

For $j=0$, we use the fact that 
\begin{equation}
 \mathds{1}_{\Omega(b) \leq \ell_j}= \frac{1}{2 \pi i} \oint_{|z|<1} \frac{z^{\Omega(b)}}{(1-z)z^{\ell_j}}  \, \frac{dz}{z}, \label{number_primes}
 \end{equation} and then
\begin{align*}
T(0) = \frac{1}{(2 \pi i)^2}& \oint \oint \\
&\sum_{\substack{ p |g_0 m_0 n_0 \Rightarrow p \in I_0 \\ (m_0,n_0)=1 }} \frac{a^2(g_0;K) a(m_0; K) a(n_0;K) \lambda(m_0) \lambda(n_0) \nu(g_0 m_0) \nu(g_0 n_0) \eta(g_0 m_0 n_0)}{g_0 m_0 n_0}  \\
& \times \frac{z_1^{\Omega(g_0 m_0)} z_2^{\Omega(g_0 n_0)}}{(1-z_1)(1-z_2) (z_1 z_2)^{\ell_0}} \, \frac{dz_1}{z_1} \frac{dz_2}{z_2},
\end{align*}
where we are integrating along circles with $|z_1|, |z_2|<1$; say $|z_1|=|z_2|=1/2$. We have that
\begin{align*}
\sum_{\substack{ p |g_0 m_0 n_0 \Rightarrow p \in I_0 \\ (m_0,n_0)=1 }} & \frac{a^2(g_0;K) a(m_0; K) a(n_0;K) \lambda(m_0) \lambda(n_0) \nu(g_0 m_0) \nu(g_0 n_0) \eta(g_0 m_0 n_0)}{g_0 m_0 n_0}  z_1^{\Omega(g_0 m_0)} z_2^{\Omega(g_0 n_0)} \\
&= \prod_{p \in I_0} \Bigg( 1 + \frac{a^2(p;K) z_1 z_2 \eta(p)}{p} - \frac{a(p;K) z_1 \eta(p)}{p} - \frac{a(p;K) z_2 \eta(p)}{p} \\
&+ \sum_{\substack{i+j+k \geq 2 \\ jk=0}} \frac{a(p^{2i+j+k};K) (-1)^{j+k} \eta(p^{i+j+k}) z_1^{i+j} z_2^{i+k}}{(i+j)! (i+k)! p^{i+j+k}} \Bigg).
\end{align*}
In the integral over $z_1$ we shift the contour of integration to $|z_1|=3$ and encounter a pole at $z_1=1$. We note that the integral over the new contour is bounded by 
\begin{align}
\frac{1}{(3/2)^{\ell_0}} & \prod_{p \in I_0} \Big(  1+ \frac{5}{p} + \sum_{b\geq 2} \frac{3^a}{p^a a!} \sum_{i \leq a} \frac{1}{i! 2^i} + \sum_{b \geq 2} \frac{1}{(2p)^b b!} \sum_{i \leq b} \frac{3^i}{i!}\Big) \nonumber \\
& \leq \frac{1}{(3/2)^{\ell_0}} \prod_{p \in I_0}  \Big(  1+ \frac{5}{p} + 2 \sum_{b\geq 2} \frac{3^b}{p^b b!}  +\frac{e^3}{4} \sum_{b \geq 2} \frac{1}{p^b b!}\Big) \nonumber \\
& \leq \frac{1}{(3/2)^{\ell_0}} \exp \Big( \sum_{p \in I_0}  \frac{5}{p}+ 2 \sum_{b \geq 2} \frac{3^b}{b!} \sum_{p \in I_0} \frac{1}{p^b}+ \frac{e^3}{4} \sum_{b \geq 2} \frac{1}{b!} \sum_{p \in I_0} \frac{1}{p^b}  \Big) \Big) \nonumber \\
& \ll \frac{1}{(3/2)^{\ell_0}} \exp \Big( 5 \log \log q^{\beta_0} + \frac{e^3 \pi^2}{3}+ \frac{e^4 \pi^2}{24} \Big) \nonumber \\
& \ll \exp \Big( - c_0 \frac{(\log q)^{3/4}}{(\log \log q)^{3/4}} \Big),
\label{et_bd}
\end{align}
for some $c_0>0$, where we used the fact that $\ell_0 \asymp (\log q)^{3/4}/(\log \log q)^{3/4}$. We then get that
\begin{align*}
T(0) &= \frac{1}{2 \pi i} \oint \prod_{p \in I_0} \Bigg( 1 + \frac{a^2(p;K)  z_2 \eta(p)}{p} - \frac{a(p;K)  \eta(p)}{p} - \frac{a(p;K) z_2 \eta(p)}{p} \\
&+ \sum_{\substack{i+j+k \geq 2 \\ jk=0}} \frac{a(p^{2i+j+k};K) (-1)^{j+k} \eta(p^{i+j+k})  z_2^{i+k}}{(i+j)! (i+k)! p^{i+j+k}} \Bigg)  \frac{ 1}{(1-z_2) z_2^{\ell_0}} \, \frac{dz_2}{z_2}\\
&+ O \Big(  \exp \Big( - c_0 \frac{(\log q)^{3/4}}{(\log \log q)^{3/4}} \Big)\Big).
\end{align*}

In the integral over $z_2$, we shift the contour of integration to $|z_2|=3$ and encounter the pole at $z_2=1$. The integral over the new contour will similarly be bounded by \eqref{et_bd}. Hence we get that 
\begin{align}
T(0) &= \prod_{p \in I_0} \Bigg( 1 + \frac{a^2(p;K) \eta(p)}{p} - \frac{2a(p;K)  \eta(p)}{p} + \sum_{\substack{i+j+k \geq 2 \\ jk=0}} \frac{a(p^{2i+j+k};K) (-1)^{j+k}\eta(p^{i+j+k}) }{(i+j)! (i+k)! p^{i+j+k}} \Bigg) \nonumber \\
 & + O \Big(  \exp \Big( - c_0 \frac{(\log q)^{3/4}}{(\log \log q)^{3/4}} \Big)\Big).\label{t0}
\end{align}
Let $U(0)$ denote the Euler product above. 
Using equations \eqref{uj}, \eqref{prod_tj}, \eqref{t0} we get that 
\begin{align}
\Big| & \prod_{j=0}^K T(j) \Big| \geq \Bigg| \prod_{p \leq q^{\beta_K}} \Bigg( 1 + \frac{a^2(p;K)  \eta(p)}{p} - \frac{2a(p;K)  \eta(p)}{p} \nonumber\\
& + \sum_{\substack{i+j+k \geq 2 \\ jk=0}} \frac{a(p^{2i+j+k};K) (-1)^{j+k} \eta(p^{i+j+k}) }{(i+j)! (i+k)! p^{i+j+k}} \Bigg) \Bigg| \nonumber \\
&  \times \Big( 1 + \frac{ O\Big( \exp \Big( - c_0 \frac{(\log q)^{3/4}}{(\log \log q)^{3/4}} \Big)\Big)} {|U(0)|}\Big)   \prod_{j=1}^K \Bigg(1 - \frac{1}{2^{\ell_j} |U(j)|} \prod_{p \in I_j} \Big(1-\frac{4}{p} \Big)^{-1} \Big( 1 - \frac{2}{p} \Big)^{-2} \Bigg). \label{et_bd2}
\end{align}

We first bound the second product in \eqref{et_bd2}. 
Using the fact that for $j \geq 1$ we have
\begin{align*}
|U(j)| & \geq \prod_{p \in I_j} \Big(1 -\frac{3}{p} \Big) \Big( 1 - \Big(1 - \frac{3}{p} \Big)^{-1} \sum_{\substack{i+j+k \geq 2 \\ jk=0}} \frac{1 }{(i+j)! (i+k)! p^{i+j+k}}\Big) \\
& \geq \prod_{p \in I_j} \Big(1 -\frac{3}{p} \Big) \Big( 1 - \frac{5}{2} \sum_{\substack{i+j+k \geq 2 \\ jk=0}} \frac{1 }{(i+j)! (i+k)! p^{i+j+k}}\Big),
\end{align*}
it follows that
\kommentar{\begin{align*}
\prod_{j=1}^K &\Bigg(1 - \frac{1}{2^{\ell_j} |U(j)|} \prod_{p \in I_j} \Big(1-\frac{4}{p} \Big)^{-1} \Big( 1 - \frac{2}{p} \Big)^{-2} \Bigg) \\
&  \geq \prod_{j=1}^K \Big(1 - \frac{1}{2^{\ell_j} }\prod_{p \in I_j} \Big(1 -\frac{3}{p} \Big)^{-1} \Big( 1 - \frac{5}{2} \sum_{\substack{i+j+k \geq 2 \\ jk=0}} \frac{1 }{(i+j)! (i+k)! p^{i+j+k}}\Big)^{-1} \Big(1-\frac{4}{p} \Big)^{-1} \Big( 1 - \frac{2}{p} \Big)^{-2} \Big) \\
& \geq \prod_{j=1}^K \Big( 1 - \frac{1}{2^{\ell_j}}\Big).
\end{align*}}
\begin{align}
&\prod_{j=1}^K\Bigg(1 - \frac{1}{2^{\ell_j} |U(j)|} \prod_{p \in I_j} \Big(1-\frac{4}{p} \Big)^{-1} \Big( 1 - \frac{2}{p} \Big)^{-2} \Bigg)  \nonumber\\
&\  \geq \prod_{j=1}^K \Big(1 - \frac{1}{2^{\ell_j} }\prod_{p \in I_j} \Big(1 -\frac{3}{p} \Big)^{-1} \Big( 1 - \frac{5}{2} \sum_{\substack{i+j+k \geq 2 \\ jk=0}} \frac{1 }{(i+j)! (i+k)! p^{i+j+k}}\Big)^{-1} \Big(1-\frac{4}{p} \Big)^{-1} \Big( 1 - \frac{2}{p} \Big)^{-2} \Big)  \nonumber \\ 
&\ \geq \prod_{j=1}^K \Big( 1 - \frac{1}{2^{\ell_j}} \sum_{p \in I_j} \Big(  \frac{22}{p} + 10 \sum_{b \geq 2} \frac{1}{p^b}\Big) \Big) \geq \prod_{j=1}^K \Big( 1 - \frac{1}{2^{\ell_j}} \sum_{p \in I_j} \Big(  \frac{22}{p} + \frac{20}{p^2} \Big) \Big)  \nonumber \\
&\ \geq \prod_{j=1}^K \Big( 1 - \frac{C}{2^{\ell_j}}\Big),  \label{prod_jbig}
\end{align} 
where $C:= 22+20\zeta(2).$
Similarly as above, we also get that
\begin{align}
 & 1 + \frac{ O\Big( \exp \Big( - c_0 \frac{(\log q)^{3/4}}{(\log \log q)^{3/4}} \Big)\Big)} {|U(0)|} \label{bound_interm}\\
& \geq 1 + O\Big( \exp \Big( - c_0 \frac{(\log q)^{3/4}}{(\log \log q)^{3/4}}\Big) \Big) \prod_{\substack{p \in I_0 }} \Big( 1 + \frac{3}{p} + \sum_{b \geq 2} \frac{1}{p^b b!} \sum_{i=0}^b \frac{1}{i!}\Big)^{-1}. \nonumber 
\end{align}
Note that
\begin{align*}
\prod_{\substack{p \in I_0 }} \Big( 1 + \frac{3}{p} + \sum_{b \geq 2} \frac{1}{p^b b!} \sum_{i=0}^b \frac{1}{i!}\Big) \leq \prod_{p \in I_0} \Big(1+ \frac{3}{p} + \frac{e}{2(p^2-p)} \Big) \leq  \prod_{p \in I_0} \Big( 1+ \frac{5}{p} \Big),
\end{align*}
so 
\begin{align*}
& 1 + \frac{ O\Big( \exp \Big( - c_0 \frac{(\log q)^{3/4}}{(\log \log q)^{3/4}} \Big)\Big)} {|U(0)|} \geq 1 + O\Big( \exp \Big( - c_0 \frac{(\log q)^{3/4}}{(\log \log q)^{3/4}}\Big) \Big) (\log q^{\beta_0})^5 \geq 1/2.
\end{align*}

Using this and \eqref{prod_jbig}, it follows that
\begin{equation*}
\Big( 1 + \frac{ O\Big( \exp \Big( - c \frac{(\log q)^{3/4}}{(\log \log q)^{3/4}} \Big)\Big)} {|U(0)|}\Big)   \prod_{j=1}^K \Bigg(1 - \frac{1}{2^{\ell_j} |U(j)|} \prod_{p \in I_j} \Big(1-\frac{4}{p} \Big)^{-1} \Big( 1 - \frac{2}{p} \Big)^{-2} \Bigg) \gg 1. 
\end{equation*}

Combining the lower bound above with \eqref{et_bd2} and \eqref{mol_mom}, we have that
\begin{align}
|\mathcal{M}| \gg  \Big|L(1,\eta) & \prod_{p \leq q^{\beta_K}} \Bigg( 1 + \frac{a^2(p;K)  \eta(p)}{p} - \frac{2a(p;K)  \eta(p)}{p}\nonumber\\
& + \sum_{\substack{i+j+k \geq 2 \\ jk=0}} \frac{a(p^{2i+j+k};K) (-1)^{j+k} \eta(p^{i+j+k}) }{(i+j)! (i+k)! p^{i+j+k}} \Bigg) \Big|. \label{intermediary} 
\end{align} 
Let
$$\mathcal{B} =L(1,\eta)  \prod_{p \leq q^{\beta_K}} \Bigg( 1 + \frac{a^2(p;K)  \eta(p)}{p} - \frac{2a(p;K)  \eta(p)}{p} + \sum_{\substack{i+j+k \geq 2 \\ jk=0}} \frac{a(p^{2i+j+k};K) (-1)^{j+k} \eta(p^{i+j+k}) }{(i+j)! (i+k)! p^{i+j+k}} \Bigg). $$
Let $E_p$ denote the term corresponding to $i+j+k \geq 2$ in the sum above.
Using \cite[Lemma $8.2$]{GS01} with $y=q^{\beta_K}$ and under GRH, we have that
$$ \log L(1,\eta) = \sum_{n \leq q^{\beta_K}} \frac{\Lambda(n) \eta(n)}{n \log n}  + O \Big(  (\log q)^3 q^{-\beta_K/2}\Big).$$
Hence it follows that
\begin{align*}
\log  &\mathcal{B} = \sum_{p \leq q^{\beta_K}} \frac{\eta(p)}{p} + \sum_{\substack{ p^b \leq q^{\beta_K} \\ b \geq 2}} \frac{\eta(p^b)}{b p^b} + \sum_{5\leq  p \leq q^{\beta_K}} \sum_{b=1}^{\infty} \frac{(-1)^{b+1}}{b}  \Big(\frac{a^2(p;K)  \eta(p)}{p} - \frac{2a(p;K)  \eta(p)}{p} +E_p \Big)^b \\
&+ \sum_{p \in \{2,3\}} \log \Big( 1 + \frac{a^2(p;K)  \eta(p)}{p} - \frac{2a(p;K)  \eta(p)}{p} + \sum_{\substack{i+j+k \geq 2 \\ jk=0}} \frac{a(p^{2i+j+k};K) (-1)^{j+k} \eta(p^{i+j+k}) }{(i+j)! (i+k)! p^{i+j+k}} \Big) \\
&= \sum_{5 \leq p \leq q^{\beta_K}} \frac{(1-a(p;K))^2 \eta(p)}{p} +\sum_{p \in \{2,3\}} \frac{\eta(p)}{p}+ \sum_{\substack{ p^b \leq q^{\beta_K} \\ b \geq 2}} \frac{\eta(p^b)}{b p^b}+ \sum_{5 \leq p \leq q^{\beta_K}} E_p \\
&+\sum_{5 \leq p \leq q^{\beta_K}} \sum_{b\geq 2} \frac{(-1)^{b+1}}{b}  \Big(\frac{a^2(p;K)  \eta(p)}{p} - \frac{2a(p;K)  \eta(p)}{p} +E_p \Big)^b \\
&+\sum_{p \in \{2,3\}} \log \Big( 1 + \frac{a^2(p;K)  \eta(p)}{p} - \frac{2a(p;K)  \eta(p)}{p} + \sum_{\substack{i+j+k \geq 2 \\ jk=0}} \frac{a(p^{2i+j+k};K) (-1)^{j+k} \eta(p^{i+j+k}) }{(i+j)! (i+k)! p^{i+j+k}} \Big) \\
&+   O \Big(  (\log q)^3 q^{-\beta_K/2}\Big).
\end{align*}
Exponentiating, we get that
\begin{align}
|\mathcal{B}| &= \Big(1  + O \Big(  (\log q)^3 q^{-\beta_K/2}\Big)  \Big)\exp \Bigg(  \Re \Big(\sum_{5 \leq p \leq q^{\beta_K}} \frac{(1-a(p;K))^2 \eta(p)}{p} + \sum_{p \in \{2,3\}} \frac{\eta(p)}{p}+\sum_{\substack{ p^b \leq q^{\beta_K} \\ b \geq 2}} \frac{\eta(p^b)}{b p^b} \nonumber \\
&+ \sum_{p \leq q^{\beta_K}} E_p +\sum_{p \leq q^{\beta_K}} \sum_{b=2}^{\infty} \frac{(-1)^{b+1}}{b}  \Big(\frac{a^2(p;K)  \eta(p)}{p} - \frac{2a(p;K)  \eta(p)}{p} +E_p \Big)^b \Big) \Bigg) \nonumber \\
& \times \prod_{p \in \{2,3\}}  \Big| 1 + \frac{a^2(p;K)  \eta(p)}{p} - \frac{2a(p;K)  \eta(p)}{p} + \sum_{\substack{i+j+k \geq 2 \\ jk=0}} \frac{a(p^{2i+j+k};K) (-1)^{j+k} \eta(p^{i+j+k}) }{(i+j)! (i+k)! p^{i+j+k}} \Big|. \label{b_exp}
\end{align}

Now we use the fact that
$$ 1-a(p;K) \leq \frac{(1+\lambda) \log p}{\beta_K \log q},$$ so
\begin{align}
\label{re1}
\Re \Big(\sum_{5 \leq p \leq q^{\beta_K}} \frac{(1-a(p;K))^2 \eta(p)}{p} \Big) \geq - \frac{(1+\lambda)^2}{\beta_K^2 (\log q)^2} \sum_{p \leq q^{\beta_K}} \frac{(\log p)^2}{p} \geq -(1+\lambda)^2.
\end{align}
We also have that 
\begin{align}
\Re \Big( \sum_{\substack{ p^b \leq q^{\beta_K} \\ b \geq 2}} \frac{\eta(p^b)}{b p^b} \Big) \geq -\sum_p \frac{1}{p^2-1},
\label{re2}
\end{align}
and
\begin{align}
\label{re3}
\Re \Big( \sum_{5 \leq p \leq q^{\beta_K}} E_p \Big) \geq -2 \sum_{p>3} \sum_{b \geq 2} \frac{1}{b! p^b} \sum_{i=0}^b \frac{1}{i!} \geq -2 \sum_{p \geq 5} \sum_{b \geq 2} \frac{b+1}{b! p^b} \geq - \sum_{p \geq 5}\frac{6}{p^2-1}.
\end{align} 
Furthermore, we have that 
\begin{align}
\Re & \Big(\sum_{5 \leq p \leq q^{\beta_K}} \sum_{b=2}^{\infty} \frac{(-1)^{b+1}}{b}  \Big(\frac{a^2(p;K)  \eta(p)}{p} - \frac{2a(p;K)  \eta(p)}{p} +E_p \Big)^b \Big) \Big) \nonumber \\
& \geq - \sum_{5 \leq  p \leq q^{\beta_K}} \sum_{b \geq 2} \frac{1}{b} \sum_{r=0}^b \binom{b}{r} \Big( \frac{6}{p^2-1} \Big)^{r} 
 \frac{3^{b-r}}{p^{b-r}} \nonumber \\
 & \geq - \sum_{5 \leq  p \leq q^{\beta_K}} \sum_{b \geq 2}  \frac{3^b}{b p^b} \Big(  1+ \frac{2p}{p^2-1}\Big)^b \geq - \frac{1}{2} \sum_{p \geq 5} \sum_{b \geq 2} \Big(\frac{9}{2p}\Big)^b \geq -\frac{9^2}{4p(2p-9)}.\label{re4}
\end{align}

Note that for $p=2,3$,
$$a(p;K) = 1+ O \Big( \frac{1}{\log q} \Big),$$
so
\begin{align}
\Big| &\frac{a^2(p;K)  \eta(p)}{p} - \frac{2a(p;K)  \eta(p)}{p} + \sum_{\substack{i+j+k \geq 2 \\ jk=0}} \frac{a(p^{2i+j+k};K) (-1)^{j+k} \eta(p^{i+j+k}) }{(i+j)! (i+k)! p^{i+j+k}} \Big| \leq \frac{1+\eps}{p} + \frac{3(1+\eps)}{4p^2} \nonumber \\
&+ \sum_{\substack{i+j+k \geq 3 \\ jk=0}} \frac{1 }{(i+j)! (i+k)! p^{i+j+k}} \leq \frac{1+\eps}{2} + \frac{3(1+\eps)}{16}+ 2 \sum_{b \geq 3} \frac{1}{b!2^b} \sum_{i=0}^b \frac{1}{i!} \nonumber  \\
& \leq \frac{1+\eps}{2} + \frac{3(1+\eps)}{16} + 2 \sum_{b \geq 3} \frac{b+1}{b! 2^b} \leq  \frac{1+\eps}{2} + \frac{3(1+\eps)}{16} + 2 \Big( \frac{1}{12}+ \frac{5}{24\cdot 16} + \frac{6}{5! 2^4} \Big) \nonumber \\
& \leq \frac{86}{96}. \nonumber 
\end{align}
It follows that
\begin{align}
\label{re5}
\Big| 1 + \frac{a^2(p;K)  \eta(p)}{p} - \frac{2a(p;K)  \eta(p)}{p} + \sum_{\substack{i+j+k \geq 2 \\ jk=0}} \frac{a(p^{2i+j+k};K) (-1)^{j+k} \eta(p^{i+j+k}) }{(i+j)! (i+k)! p^{i+j+k}} \Big| \geq \frac{5}{48}.
\end{align}

Now using equations \eqref{b_exp}, \eqref{re1}, \eqref{re2}, \eqref{re3}, \eqref{re4}, \eqref{re5}, it follows that 
\begin{align*}
|\mathcal{B} | & \geq \Big(1  + O \Big(  (\log q)^3 q^{-\beta_K/2}\Big)  \Big) \Big(  \frac{5}{48}\Big)^2 \exp \Big(  - \frac{5}{6}- (1+\lambda)^2 - 7 \sum_p \frac{1}{p^2-1} -\sum_{p \geq 5}  \frac{9^2}{4p(2p-9)}\Big) \\
& \geq \Big(1  + O \Big(  (\log q)^3 q^{-\beta_K/2}\Big)  \Big) \Big(  \frac{5}{48}\Big)^2 \exp \Big(  -\frac{5}{6}- (1+\lambda)^2  - (14+ \tfrac{3}{2} \cdot 9^2) \sum_p \frac{1}{p^2} \Big) \\
& \gg 1.
\end{align*}
Using the above and equation \eqref{intermediary}, we get that
\begin{equation}
|\mathcal{M}| \gg 1.
\end{equation}
This completes the proof of the theorem.

\section{Proof of Theorem \ref{thm:nonvanishing-full-orbit} }

H\"{o}lder's inequality implies that
\begin{align}\label{HI}
&\bigg(\sum_{\substack{\chi\in\co\\L(1/2,\chi\eta_1)L(1/2,\chi\eta_2)\ne0}}1\bigg)^2\bigg(\sum_{\substack{\chi\in\co}}|L(\tfrac12,\chi\eta_1)M(\chi\eta_1)|^4\bigg)\bigg(\sum_{\substack{\chi\in\co}}|L(\tfrac12,\chi\eta_2)M(\chi\eta_2)|^4\bigg)\nonumber\\
&\qquad\qquad\geq \bigg(\sum_{\substack{\chi\in\co}}L(\tfrac12,\chi\eta_1)L(\tfrac12,\overline{\chi\eta_2})M(\chi\eta_1)M(\overline{\chi\eta_2})\bigg)^4.
\end{align}
By positivity,
\begin{align*}  
\sum_{\substack{\chi\in\co}}|L(\tfrac12,\chi\eta_1)M(\chi\eta_1)|^4\leq \sumstar_{\chi(\text{mod}\ q)}|L(\tfrac12,\chi)M(\chi)|^4.
\end{align*} 
Using Theorem \ref{mol_ub} with $v=4$, we have that
$$\sumstar_{\chi(\text{mod}\ q)}|L(\tfrac12,\chi)M(\chi)|^4 \ll q, $$ 
and since $|\co|\asymp q$, we obtain from \eqref{HI} that
\begin{equation}
\label{lbmol}
\frac{1}{|\co|}\sum_{\substack{\chi\in\co\\L(1/2,\chi\eta_1)L(1/2,\chi\eta_2)\ne0}}1\geq \bigg(\frac{1}{|\co|}\,\sum_{\substack{\chi\in\co}}L(\tfrac12,\chi\eta_1)L(\tfrac12,\overline{\chi\eta_2})M(\chi\eta_1)M(\overline{\chi\eta_2})\bigg)^2.
\end{equation}
Using Theorem \ref{mol_lb} finishes the proof.
\qed

\subsection*{Acknowledgements} HMB was supported by the Dame Kathleen Ollerenshaw (DKO) Mathematics Travel Fund of the University of Manchester when the work commenced. AF was supported by the National Science Foundation (DMS-2101769 and CAREER grant DMS-2339274). Part of this work was done while HTN was visiting the Vietnam Institute for Advanced Study in Mathematics (VIASM); HTN would like to thank VIASM for their hospitality and partial support.

\bibliographystyle{amsplain}

\bibliography{galois-orbits-bib}

\providecommand{\bysame}{\leavevmode\hbox to3em{\hrulefill}\thinspace}
\providecommand{\MR}{\relax\ifhmode\unskip\space\fi MR }
\providecommand{\MRhref}[2]{%
  \href{http://www.ams.org/mathscinet-getitem?mr=#1}{#2}
}
\providecommand{\href}[2]{#2}
\begin{thebibliography}{10}

\bibitem{Avc25}
O.~Avci, \emph{Torsion of rational elliptic curves over the
  {$\mathbb{Z}_p$}-extensions of quadratic fields}, preprint,
  {arXiv}:2505.04149.

\bibitem{BY10}
S.~Baier and M.~P. Young, \emph{Mean values with cubic characters}, J. Number
  Theory \textbf{130} (2010), 879--903.

\bibitem{Bui12}
H.~M. Bui, \emph{Non-vanishing of {D}irichlet {$L$}-functions at the central
  point}, Int. J. Number Theory \textbf{8} (2012), 1855--1881.

\bibitem{BELP25}
H.~M. Bui, N.~Evans, S.~Lester, and K.~Pratt, \emph{Weighted central limit
  theorems for central values of {{\(L\)}}-functions}, J. Eur. Math. Soc.
  \textbf{27} (2025), 2477--2529.

\bibitem{Cha09}
V.~Chandee, \emph{Explicit upper bounds for {{\(L\)}}-functions on the critical
  line}, Proc. Amer. Math. Soc. \textbf{137} (2009), 4049--4063.

\bibitem{Chi02}
G.~Chinta, \emph{Analytic ranks of elliptic curves over cyclotomic fields}, J.
  Reine Angew. Math. \textbf{544} (2002), 13--24.

\bibitem{Cho65}
S.~Chowla, \emph{The {Riemann} hypothesis and {Hilbert}'s tenth problem},
  London and {Glasgow}: {Blackie} \& {Son} {Ltd}. {XV}, 119 pages, 1965.

\bibitem{DDDS24}
C.~David, A.~de~Faveri, A.~Dunn, and J.~Stucky, \emph{Non-vanishing for cubic
  {Hecke} ${L}$-functions}, preprint, {arXiv}:2410.03048.

\bibitem{DFL21}
C.~David, A.~Florea, and M.~Lalin, \emph{Nonvanishing for cubic
  {{\(L\)}}-functions}, Forum Math. Sigma \textbf{9} (2021), 58, e69.

\bibitem{DG22}
C.~David and A.~M. G{\"u}lo{\u{g}}lu, \emph{One-level density and non-vanishing
  for cubic {{\(L\)}}-functions over the {Eisenstein} field}, Int. Math. Res.
  Not. (2022), no.~23, 18833--18873.

\bibitem{GS01}
A.~Granville and K.~Soundararajan, \emph{Large character sums}, J. Amer. Math.
  Soc. \textbf{14} (2001), 365--397.

\bibitem{Gre83}
R.~Greenberg, \emph{On the {Birch} and {Swinnerton}-{Dyer} conjecture}, Invent.
  Math. \textbf{72} (1983), 241--265.

\bibitem{Gre85}
\bysame, \emph{On the critical values of {H}ecke {$L$}-functions for imaginary
  quadratic fields}, Invent. Math. \textbf{79} (1985), 79--94.

\bibitem{Gre87}
\bysame, \emph{Non-vanishing of certain values of {{\(L\)}}-functions},
  Analytic {N}umber {T}heory and {D}iophantine {P}roblems, {Proc}. {Conf}.,
  {Stillwater}/{Okla}. 1984, {Prog}. {Math}. 70, 223-235, 1987.

\bibitem{Gre01}
\bysame, \emph{Introduction to {Iwasawa} theory for elliptic curves},
  Arithmetic {A}lgebraic {G}eometry, Providence, RI: American Mathematical
  Society (AMS), 2001, pp.~407--464.

\bibitem{Har13}
A.~J. Harper, \emph{Sharp conditional bounds for moments of the {Riemann} zeta
  function}, preprint, {arXiv}:1305.4618.

\bibitem{IS99}
H.~Iwaniec and P.~Sarnak, \emph{Dirichlet {$L$}-functions at the central
  point}, Number theory in progress, {V}ol. 2 ({Z}akopane-{K}o\'scielisko,
  1997), de Gruyter, Berlin, 1999, pp.~941--952.

\bibitem{KMN16}
R.~Khan, D.~Mili{\'c}evi{\'c}, and H.~T. Ngo, \emph{Non-vanishing of
  {Dirichlet} {$L$}-functions in {Galois} orbits}, Int. Math. Res. Not. (2016),
  no.~22, 6955--6978.

\bibitem{KMN22}
\bysame, \emph{Nonvanishing of {Dirichlet} {{\(L\)}}-functions. {II}}, Math. Z.
  \textbf{300} (2022), 1603--1613.

\bibitem{KN16}
R.~Khan and H.~T. Ngo, \emph{Nonvanishing of {Dirichlet} {{\(L\)}}-functions},
  Algebra Number Theory \textbf{10} (2016), 2081--2091.

\bibitem{LR21}
S.~Lester and M.~Radziwi{\l}{\l}, \emph{Signs of {Fourier} coefficients of
  half-integral weight modular forms}, Math. Ann. \textbf{379} (2021), no.~3-4,
  1553--1604.

\bibitem{Maz72}
B.~Mazur, \emph{Rational points of abelian varieties with values in towers of
  number fields}, Invent. Math. \textbf{18} (1972), 183--266.

\bibitem{MV02}
P.~Michel and J.~Vanderkam, \emph{Simultaneous nonvanishing of twists of
  automorphic {{\(L\)}}-functions}, Compos. Math. \textbf{134} (2002),
  135--191.

\bibitem{MR82}
H.~L. Montgomery and D.~E. Rohrlich, \emph{On the {$L$}-functions of canonical
  {Hecke} characters of imaginary quadratic fields. {II}}, Duke Math. J.
  \textbf{49} (1982), 937--942.

\bibitem{MV07}
H.~L. Montgomery and R.~C. Vaughan, \emph{Multiplicative number theory. {I}.
  {Classical} theory}, Camb. Stud. Adv. Math., vol.~97, Cambridge: Cambridge
  University Press, 2007.

\bibitem{MM23}
M.~Munsch and I.~E. Shparlinski, \emph{Moments and non-vanishing of
  {$L$}-functions over thin subgroups}, preprint, {arXiv}:2309.10207.

\bibitem{OS99}
A.~E. {\"O}zl{\"u}k and C.~Snyder, \emph{On the distribution of the nontrivial
  zeros of quadratic {{\(L\)}}-functions close to the real axis}, Acta Arith.
  \textbf{91} (1999), 209--228.

\bibitem{PY23}
I.~Petrow and M.~P. Young, \emph{The fourth moment of {Dirichlet}
  {{\(L\)}}-functions along a coset and the {Weyl} bound}, Duke Math. J.
  \textbf{172} (2023), 1879--1960.

\bibitem{QW25}
X.~Qin and X.~Wu, \emph{Non-vanishing of {Dirichlet} {$L$}-functions at the
  central point}, preprint, {arXiv}:2504.11916.

\bibitem{Rid58}
D.~Ridout, \emph{The {$p$}-adic generalization of the {T}hue-{S}iegel-{R}oth
  theorem}, Mathematika \textbf{5} (1958), 40--48.

\bibitem{Roh80b}
D.~E. Rohrlich, \emph{Galois conjugacy of unramified twists of {Hecke}
  characters}, Duke Math. J. \textbf{47} (1980), 695--703.

\bibitem{Roh80a}
\bysame, \emph{The non-vanishing of certain {Hecke} {$L$}-functions at the
  center of the critical strip}, Duke Math. J. \textbf{47} (1980), 223--232.

\bibitem{Roh80c}
\bysame, \emph{On the {$L$}-functions of canonical {Hecke} characters of
  imaginary quadratic fields}, Duke Math. J. \textbf{47} (1980), 547--557.

\bibitem{Roh84a}
\bysame, \emph{On {$L$}-functions of elliptic curves and anticyclotomic
  towers}, Invent. Math. \textbf{75} (1984), 383--408.

\bibitem{Roh84b}
\bysame, \emph{On {$L$}-functions of elliptic curves and cyclotomic towers},
  Invent. Math. \textbf{75} (1984), 409--423.

\bibitem{Roh88}
\bysame, \emph{L-functions and division towers}, Math. Ann. \textbf{281}
  (1988), 611--632.

\bibitem{Rot55}
K.~F. Roth, \emph{Rational approximations to algebraic numbers}, Mathematika
  \textbf{2} (1955), 1--20.

\bibitem{Shi76}
G.~Shimura, \emph{The special values of the zeta functions associated with cusp
  forms}, Comm. Pure Appl. Math. \textbf{29} (1976), 783--804.

\bibitem{Shi77}
\bysame, \emph{On the periods of modular forms}, Math. Ann. \textbf{229}
  (1977), 211--221.

\bibitem{Sou00}
K.~Soundararajan, \emph{Nonvanishing of quadratic {D}irichlet {$L$}-functions
  at {$s=\frac12$}}, Ann. of Math. (2) \textbf{152} (2000), 447--488.

\bibitem{Sou09}
\bysame, \emph{Moments of the {Riemann} zeta function}, Ann. of Math. (2)
  \textbf{170} (2009), 981--993.

\bibitem{Sza24}
Barnab{\'a}s Szab{\'o}, \emph{High moments of theta functions and character
  sums}, Mathematika \textbf{70} (2024), no.~2, 37 pp.

\bibitem{Zac19}
R.~Zacharias, \emph{Simultaneous non-vanishing for {Dirichlet}
  {{\(L\)}}-functions}, Ann. Inst. Fourier \textbf{69} (2019), 1459--1524.

\end{thebibliography}

\end{document}